\documentclass[11pt]{article}

\usepackage[a4paper,margin=1in]{geometry}
\usepackage{amsmath,amssymb,amsfonts,amsthm,mathtools}
\usepackage{microtype}
\usepackage{enumitem}
\usepackage{xcolor}
\usepackage{authblk}
\usepackage{hyperref}
\hypersetup{colorlinks=true,linkcolor=blue,citecolor=purple,urlcolor=teal}
\usepackage{comment}

\DeclareMathOperator{\GL}{GL}
\DeclareMathOperator{\PG}{PG}

\newcommand{\wt}{w}

\newcommand{\F}{\mathbb F}
\newcommand{\fq}{\mathbb F_q}
\newcommand{\fqm}{\mathbb F_{q^m}}
\newcommand{\fqmk}{\mathbb F_{q^m}^{k}}
\newcommand{\dimq}{\dim_{\mathbb F_q}}
\newcommand{\dimqm}{\dim_{\mathbb F_{q^m}}}
\newcommand{\Uperpp}{U^{\perp'}}
\newcommand{\Luperpp}{L_{U^{\perp'}}}
\newcommand{\Us}{U^*}
\newcommand{\Usqm}{\langle U^* \rangle_{\fqm} }
\newcommand{\Upqm}{\langle U' \rangle_{\fqm} }
\newcommand{\Sp}{S^\perp}

\newcommand{\cC}{\mathcal C}
\newcommand{\cD}{\mathcal D}

\newcommand{\Tr}{\operatorname{Tr}}
\newcommand{\Exp}{\operatorname{Exp}}

\newtheorem{theorem}{Theorem}[section]
\newtheorem{proposition}[theorem]{Proposition}
\newtheorem{lemma}[theorem]{Lemma}
\newtheorem{corollary}[theorem]{Corollary}

\theoremstyle{definition}
\newtheorem{definition}[theorem]{Definition}
\newtheorem{example}[theorem]{Example}
\theoremstyle{remark}
\newtheorem{remark}[theorem]{Remark}

\title{Few-Weight Linear Sets and Divisibility of Rank-Metric Codes}

\author[1]{Jan De Beule}
\author[1]{Leen Demuys}
\author[2]{Olga Polverino}
\author[2]{Ferdinando Zullo}

\affil[1]{Department of Mathematics and Data Science,
Vrije Universiteit Brussel, Pleinlaan 2, 1050 Brussel, Belgium}

\affil[2]{Dipartimento di Matematica e Fisica,
Universit\`a degli Studi della Campania ``Luigi Vanvitelli'',
Viale Lincoln 5, I--81100 Caserta, Italy}

\date{}

\begin{document}
\maketitle

\begin{abstract}
We study $\mathbb F_q$-linear sets whose positive point weights belong to a prescribed finite set $A$.  The integer $e=\gcd(A)$ always divides the rank of the defining subspace, but divisibility of the weights does not by itself imply linearity over $\mathbb F_{q^e}$.  We show that it does in a broad parameter range.  Let $L_U$ be an $\mathbb F_q$-linear set of rank $N$ spanning $\PG(h-1,q^m)$, $h\ge2$, without being equal to it, and assume that $q\ge m$ and that {the minimum weight of its points is $a\ge2$}.  Combining the theorem of Csajb\'ok, Marino and Pepe on the maximum field of linearity with the lower bound of Adriaensen and Santonastaso for the size of a linear set containing a point of weight one, we prove that $L_U$ is also defined by an $\mathbb F_{q^t}$-subspace, for some divisor $t$ of $m$, whose $\mathbb F_q$-dimension is exactly $N+t-a$.  Under a natural gap condition excluding divisors of $m$ in the interval $(e,N-e]$, this forces $a=e$ and the defining subspace to be $\mathbb F_{q^e}$-linear.

Through geometric duality, the point weights of the dual of a system are the deficiencies $m-w$ of the rank weights $w$ of the associated code.  For rank-metric codes in which the subcode generated by the codewords of nonmaximum rank contains codewords of maximum rank, we show that deleting suitable coordinates yields the $\mathbb F_q$-expansion of a code over $\mathbb F_{q^t}$, in which every codeword of rank $w$ becomes a codeword of rank $t\lfloor w/t\rfloor$.  This also yields congruence restrictions on the weight spectrum. We show that the codes in which the codewords of nonmaximum rank generate a subcode without codewords of maximum rank are exactly the direct sums of a simplex code and a code of length smaller than $m$.  We derive field-reduction results for two-weight and three-weight rank-metric codes, and show that no gap condition is needed when $m$ divides the length or when the code has only two nonzero weights.  Finally, we study two-dimensional codes of length at most $m$ and provide families of three-weight examples.
\end{abstract}


\section{Introduction}

Restrictions on the weights of a linear code can force strong algebraic and geometric structure. This principle is particularly visible for codes with only a few nonzero weights. In the Hamming metric, two-weight codes connect coding theory with point sets having two hyperplane intersection numbers and with strongly regular graphs; see the survey of Calderbank and Kantor~\cite{CalderbankKantor1986}. Divisibility of weights provides another source of structural constraints, with a substantial theory of its own~\cite{Ward2001}. These connections motivate the study of how much of a code's structure can be recovered from its weight spectrum and the arithmetic relations among its weights.

The rank metric offers a natural setting for this problem. Here the weight of a vector over a finite extension field is the dimension, over the base field, of the space generated by its coordinates. The foundational work of Delsarte and Gabidulin~\cite{Delsarte1978,Gabidulin1985} established an algebraic theory of codes for this metric, which also plays a central role in error control for random network coding~\cite{SilvaKschischangKoetter2008}. Since rank weight measures linear dependence over a chosen base field, restrictions on rank weights are closely connected with the interaction between different field structures.

Finite geometry makes this connection explicit. A linear code in the rank metric can be represented by a vector subspace over the base field, and codeword weights are determined by its intersections with hyperplanes~\cite{Randrianarisoa2020,AlfaranoEtAl2022}. The associated projective point sets are called \emph{linear sets}; their point weights record the dimensions of the intersections with the corresponding one-dimensional subspaces. Linear sets have been extensively studied in finite geometry~\cite{Polverino2010}, and their duality theory provides a bridge between geometric intersection properties and rank-weight distributions~\cite{SheekeyVandeVoorde2020,BorelloPolverinoZullo2025}.

This geometric approach has been particularly effective for codes with few weights. Constant-weight rank-metric codes were classified in~\cite{Randrianarisoa2020}, and two-weight rank-metric codes were described through scattered subspaces in~\cite{ZulloEtAl2024}. Three-weight families associated with clubs were studied in~\cite{MannaertSantonastasoZullo2025} (see \cite{SmaldoreZanellaZullo2026,Zullo2023} for some extensions). These results suggest investigating which features of such descriptions persist when the entire weight spectrum is restricted, without prescribing a particular geometric configuration.

A central issue is whether divisibility of weights reveals an underlying structure over a larger field. Field expansion gives a natural source of divisible rank-metric codes: expanding a code from a larger base field multiplies all rank weights by the degree of the field extension. The converse is more delicate. In~\cite{PolverinoEtAl2023} the authors identified cases in which divisibility does imply such an origin, but also constructed counterexamples. Thus numerical divisibility alone does not determine the field structure. The problem is to find additional conditions, expressed in terms of the weight spectrum and the code parameters, that make this implication valid.

In this paper we approach this problem through linear sets with restricted point weights. If the defining subspace is linear over a larger subfield, then every point weight is divisible by the degree of that subfield over the base field. We investigate when the converse holds and, when it does not hold directly, how far the defining subspace is from one with a larger field of linearity. This distinction is essential: a point set may admit a representation over a larger field even when its original defining subspace is not linear over that field.

Our first main result quantifies this distinction. Let $U$ be an $\fq$-subspace of rank $N$ defining a linear set $L_U$ that spans $\PG(h-1,q^m)$, with $h\ge2$, and is a proper subset of that space. Assume that $q\ge m$ and that the minimum positive point weight is $a\ge2$. Building on the theorem of Csajb\'ok, Marino and Pepe~\cite{CsajbokMarinoPepe2024}, we show that there is a divisor $t$ of $m$ such that $W=\langle U\rangle_{\mathbb F_{q^t}}$ defines the same point set, $L_W$ admits a point of weight one over $\F_{q^t}$ and
\[
  \dimq W=N+t-a.
\]
Moreover, $a\le t\le N-a$ (Theorem~\ref{thm:larger-field-representation}). The exact dimension follows by comparing an elementary counting bound with the lower bound of Adriaensen and Santonastaso~\cite{AdriaensenSantonastaso2024} for linear sets containing a point of weight one; see also~\cite{DeBeuleVandeVoorde2019}.

This dimension formula yields a structural criterion. Let $e$ be the greatest common divisor of the positive point weights. An elementary counting argument gives $e\mid N$. Under the preceding assumptions, if $m$ has no divisor in the interval $(e,N-e]$, then $t=a=e$ and $W=U$ (Theorem~\ref{thm:gcd-rigidity}). Consequently, the original subspace is $\mathbb F_{q^e}$-linear, and its minimum positive point weight equals the gcd of all its point weights. The arithmetic restriction therefore identifies the larger field and rules out the need to enlarge the defining subspace.

{The first application concerns two-dimensional codes of length $n$ with a codeword of weight $n$.  The hyperplanes of $\PG(1,q^m)$ are points, so by
\eqref{eq:code-system} the positive point weights of an associated system $U$ are
the differences $n-d_i$ between the length and the nonzero weights $d_i<n$ of the
code, and the minimum one is the gap $a=n-d_r$ between the two largest weights.
A nontrivial $\fq$-expansion of such a code always has $a\ge2$.  Conversely, if
$q\ge m$ and $a\ge2$, applying Theorem~\ref{thm:larger-field-representation} to $U$
shows that the code is obtained by deleting $t-a$ coordinates from the
$\fq$-expansion of a $[d_r/t+1,2]_{q^m/q^t}$ code, for some divisor $t$ of $m$ with
$a\le t\le\min\{d_r,m-d_r\}$ and $t\mid d_r$; a codeword of rank $w$ arises from a
codeword of rank $t\lceil w/t\rceil$ (Theorem~\ref{thm:antipodal-puncturing}).  This
gives congruence restrictions on the weights, the inequality $d_1+d_r\ge n$ and,
for three-weight codes with $t>a$, an explicit weight distribution
(Corollary~\ref{cor:antipodal-puncturing-three}).  We complement these results
with several families of three-weight $[n,2]_{q^m/q}$ codes, which exist for every
prime power $q$.  They show that, when $a=1$, no arithmetic restriction survives,
even if the gap between the two smallest weights is large, and that the case $t>a$
really occurs: Family~$\mathcal F_4$ provides codes whose first gap $t-1$ can be as
large as $m/2-1$ and which are not expansions from any larger field.  This is in
sharp contrast with the two-weight case, in which $q\ge m$ and $a\ge2$ force the
code itself to be an expansion (Theorem~\ref{thm:two-weight-short}).}

For rank-metric codes of larger dimension, the link with these results is provided by geometric duality: the positive point weights of the dual of an associated system are precisely the deficiencies $m-w$ of the codeword ranks $w<m$. For a nondegenerate code of length at least $m$, a useful dichotomy is determined by the subcode generated by its codewords of rank smaller than $m$. If this subcode contains no codeword of rank $m$, we characterize the original code, up to equivalence, as the direct sum of a simplex code and a code of length smaller than $m$ (Theorem~\ref{thm:simplex-short-code} and Corollary~\ref{cor:direct-sum-weights}). This reduces the determination of its weight spectrum to the study of the weights of a shorter code.

In the complementary case, that is when the subcode generated by its codewords of rank smaller than $m$ contains a codeword of rank $m$, assume $q\ge m$ and let $a\ge2$ be the difference between $m$ and the largest nonmaximum rank weight. We prove that, after a rank-metric equivalence, deleting $t-a$ suitable coordinates yields the $\fq$-expansion of a code whose rank is measured over $\mathbb F_{q^t}$, for some divisor $t$ of $m$. The effect on every codeword is explicit: a word of rank $w$ becomes a word of rank $t\lfloor w/t\rfloor$ (Theorem~\ref{thm:larger-field-puncturing}). This gives both congruence restrictions on the original weight spectrum and the weight distribution after puncturing. Under the corresponding divisor-gap condition, no coordinates need to be deleted, and the original code itself arises by field expansion.

We apply these results to codes with two or three nonzero weights. The two-weight case admits stronger conclusions: in the field-size range $q\ge m$, the relevant field-expansion results require no divisor-gap assumption, with the direct-sum case treated separately. A further simplification occurs when $m$ divides the code length: the expansion theorem then holds without a divisor-gap condition or a restriction on $q$ (Theorem~\ref{thm:length-multiple-m}). 
{In Section~\ref{sec:few-weight} we apply these results to codes with two or three
nonzero weights.  For two-weight codes no divisor-gap condition is needed: when
$q\ge m$, a linear set whose points all have the same weight $a\ge2$ is forced to
be $\mathbb F_{q^a}$-linear (Lemma~\ref{lem:one-weight-linear-set}).  If the length
is at least $m$, then either $h(\cC)=1$ and $\cC$ is the direct sum of a simplex
code and a one-dimensional code, or $h(\cC)\ge2$ and, when $q\ge m$ or $m$ divides
the length, $\cC$ is the $\fq$-expansion of a code over $\mathbb F_{q^a}$, where
$a$ is the difference between the two weights (Theorem~\ref{thm:two-weight}).  If
the length is at most $m$, the code is two-dimensional and, when $q\ge m$ or the
length equals $m$, it is again an expansion (Theorem~\ref{thm:two-weight-short}).
This recovers the results of~\cite{ZulloEtAl2024}, shows that the description given
there for length larger than $m$ requires the additional assumption $h(\cC)\ge2$,
and removes the restriction on $q$ when $m$ divides the length.  For three-weight
codes of length at least $m$, the direct-sum case reduces to two-weight codes of
length smaller than $m$, and its weight distribution is determined by the
parameters (Theorem~\ref{thm:three-weight-direct}); in the remaining case we obtain
the puncturing and rigidity statements described above
(Theorem~\ref{thm:three-weight-not-direct}).  Together with the families of
Section~\ref{sec:antipodal}, these results show that three-weight codes are much
less rigid than two-weight codes.
Our results thus identify parameter regimes in which the question raised
in~\cite{PolverinoEtAl2023}, namely whether an $e$-divisible code arises from a
code over $\mathbb F_{q^e}$, has a positive answer, under hypotheses stated
directly in terms of the weight spectrum and the divisors of $m$.}

\paragraph{Organization.}
Section~\ref{sec:preliminaries} collects the preliminaries on linear sets, rank-metric codes, $q$-systems and field expansion.  In Section~\ref{sec:linear-sets} we prove Theorems~\ref{thm:larger-field-representation} and~\ref{thm:gcd-rigidity} on the fields of linearity of linear sets with restricted point weights.  Section~\ref{sec:antipodal} deals with two-dimensional codes of length at most $m$ and contains several families of three-weight examples.  Section~\ref{sec:codes} translates the results of Section~\ref{sec:linear-sets} to rank-metric codes, distinguishing between direct sums of a simplex code and a short code (Section~\ref{subsec:direct-sum}) and the remaining codes (Section~\ref{subsec:not-direct-sum}).  Finally, Section~\ref{sec:few-weight} collects the results on codes with two or three nonzero weights.

\section{Preliminaries}\label{sec:preliminaries}

Throughout the paper, $q$ is a prime power, $m$ and $k$ are positive integers, and $V$ is a $k$-dimensional vector space over $\fqm$.  We identify points of $\PG(V,\fqm)=\PG(k-1,q^m)$ with one-dimensional $\fqm$-subspaces of $V$.

\subsection{Linear sets}

We begin by recalling the basic terminology and notation for $\mathbb F_q$-linear sets that will be used throughout the paper.

\begin{definition}
Let $U$ be an $\fq$-subspace of $V$.  The $\fq$-\textbf{linear set} defined by $U$ is
\[
  L_U=\{\langle u\rangle_{\fqm}:u\in U\setminus\{0\}\}
  \subseteq \PG(V,\fqm).
\]
Its \textbf{rank} is $\dimq U$.  For a point $P=\langle v\rangle_{\fqm}$, its \textbf{weight} with respect to $U$ is
\[
  w_U(P)=\dimq(U\cap \langle v\rangle_{\fqm}).
\]
We set $w_U(P)=0$ when $P\notin L_U$.
\end{definition}

We refer to \cite{Polverino2010} for background on linear sets.  The following vector-counting identity will be used repeatedly.  If $N_i$ denotes the number of points of weight $i$, then
\begin{equation}\label{eq:vector-counting}
  q^{\dimq U}-1=\sum_{i\ge1}N_i(q^i-1).
\end{equation}
In particular, if every positive point weight is at least $a$, then
\begin{equation}\label{eq:size-upper-minweight}
  |L_U|\le \frac{q^{\dimq U}-1}{q^a-1}.
\end{equation}
We shall usually replace $V$ with the $\fqm$-span of $U$.  Thus, when convenient, we assume that
\[
  \langle U\rangle_{\fqm}=V.
\]
If $h=\dimqm\langle U\rangle_{\fqm}$, then $L_U$ spans a copy of $\PG(h-1,q^m)$.

The next elementary lemma provides a characterization of when a linear set fills its span.

\begin{lemma}\label{lem:covering}
Let $S$ be an $h$-dimensional $\fqm$-vector space, $h\ge1$, and let $U\le_{\fq}S$.  Then $L_U=\PG(S,\fqm)$ if and only if $\dimq U>(h-1)m$.
\end{lemma}

\begin{proof}
If $\dimq U>(h-1)m$, then, for every one-dimensional $\fqm$-subspace $P\le S$, Grassmann's formula gives
\[
  \dimq(U\cap P)\ge \dimq U+m-hm>0,
\]
so $P\in L_U$.  Conversely, if $L_U=\PG(S,\fqm)$, then every point has weight at least one, and \eqref{eq:vector-counting} gives
\[
  q^{\dimq U}-1
  \ge \frac{q^{hm}-1}{q^m-1}(q-1)
  =(q^{(h-1)m}+\cdots+q^m+1)(q-1)
  \ge q^{(h-1)m}
  >q^{(h-1)m}-1.
\]
Hence $\dimq U>(h-1)m$.
\end{proof}

Fix a nondegenerate symmetric $\fqm$-bilinear form
\[
  \sigma:V\times V\longrightarrow \fqm
\]
and define the nondegenerate $\fq$-bilinear form
\[
  \sigma'(x,y)=\Tr_{q^m/q}(\sigma(x,y)).
\]
We write $\perp$ for orthogonality with respect to $\sigma$ and $\perp'$ for orthogonality with respect to $\sigma'$.  Recall that $X^{\perp'}=X^\perp$ for every $\fqm$-subspace $X$ of $V$; see \cite{Polverino2010}.  When $V=\fqm^k$, we take $\sigma$ to be the standard inner product $\sigma(x,y)=\sum_i x_iy_i$.

\begin{proposition}\cite[Property 2.6]{Polverino2010}\label{prop:weight-duality}
Let $U\le_{\fq}V$ have rank $n$, and let $X\le_{\fqm}V$ have $\fqm$-dimension $s$.  Then
\begin{equation}\label{eq:general-duality}
\dimq(U^{\perp'}\cap X^\perp)
  =\dimq(U\cap X)+km-n-sm.
\end{equation}
In particular, for every $\fqm$-hyperplane $H$,
\begin{equation}\label{eq:hyperplane-point-duality}
  w_{U^{\perp'}}(H^\perp)
  =w_U(H)+m-n.
\end{equation}
\end{proposition}

\begin{remark}\label{rem:baseline}
For every $\fqm$-hyperplane $H$, Grassmann's formula yields
\[
  \dimq(U\cap H)\ge n-m,
\]
which is informative as soon as $n>m$.  Formula \eqref{eq:hyperplane-point-duality} shows that the point weights of $U^{\perp'}$ are exactly the excesses of the hyperplane weights of $U$ above this baseline.
\end{remark}

We shall repeatedly use that trace duality preserves linearity over subfields.

\begin{lemma}\label{lem:perp-linear}
Let $s\mid m$ and let $X$ be an $\mathbb F_{q^s}$-subspace of $V$.  Then $X^{\perp'}$ is also an $\mathbb F_{q^s}$-subspace of $V$.
\end{lemma}

\begin{proof}
Let $u\in X^{\perp'}$, $\lambda\in\mathbb F_{q^s}$ and $x\in X$.  Since $\lambda x\in X$ and $\sigma$ is $\fqm$-bilinear,
\[
  \sigma'(\lambda u,x)
  =\Tr_{q^m/q}(\lambda\,\sigma(u,x))
  =\Tr_{q^m/q}(\sigma(u,\lambda x))
  =\sigma'(u,\lambda x)=0.
\]
Thus $\lambda u\in X^{\perp'}$.
\end{proof}

\subsection{Rank-metric codes and $q$-systems}

We briefly recall the correspondence between linear rank-metric codes
and $q$-systems, which allows us to translate rank weights of codewords
into intersection dimensions with hyperplanes.  This correspondence
will be used repeatedly in the sequel.

For a vector $c=(c_1,\ldots,c_n)\in \fqm^n$,
its \textbf{rank weight} is
\[
\wt_q(c)
=
\dim_{\fq}\langle c_1,\ldots,c_n\rangle_{\fq}.
\]
If $q$ is clear from the context, we will just write $\wt(c)$.
The \textbf{rank distance} on $\fqm^n$ is defined by $d_{\rm R}(c,c')=\wt(c-c')$.
An $[n,k]_{q^m/q}$ \textbf{rank-metric code} is an
$\fqm$-subspace $\cC\le \fqm^n$ of dimension $k$.  Its minimum rank distance is
\[
d(\cC)
=
\min\{\wt(c):c\in\cC\setminus\{0\}\}.
\]
Clearly, every codeword has rank weight at most $\min\{n,m\}$.

Let $G\in\fqm^{k\times n}$ be a generator matrix of $\cC$, and
write its columns as
\[
G=(u_1\ \cdots\ u_n),
\qquad
u_i\in\fqm^k.
\]
We say that $\cC$ is \textbf{nondegenerate} if the columns of $G$ are
linearly independent over $\fq$.  This condition does not depend on
the choice of the generator matrix.  Equivalently, if $U=\langle u_1,\ldots,u_n\rangle_{\fq}\le\fqm^k$,
then $\dimq U=n$.
Since $G$ has rank $k$ over $\fqm$, one also has $\langle U\rangle_{\fqm}=\fqm^k$.

An $\fq$-subspace $U\le\fqm^k$ satisfying
$\dimq U=n$ and $\langle U\rangle_{\fqm}=\fqm^k$ is called an $[n,k]_{q^m/q}$ \textbf{$q$-system}.  Thus every nondegenerate rank-metric code determines a $q$-system by taking the
$\fq$-span of the columns of a generator matrix.  Conversely, given
an $[n,k]_{q^m/q}$ system $U$, choosing an $\fq$-basis
$u_1,\ldots,u_n$ of $U$ as the columns of a matrix $G$ produces a
nondegenerate $[n,k]_{q^m/q}$ rank-metric code.

Changing the $\fq$-basis of $U$ replaces $G$ by $GA$, for some $A\in\GL(n,q)$.
Accordingly, throughout the paper we regard codes related by such a
transformation as equivalent.  On the geometric side, changing the
generator matrix of the same code by left multiplication with an
element of $\GL(k,q^m)$ replaces $U$ by a projectively equivalent
$q$-system.  In this sense, the correspondence between nondegenerate
rank-metric codes and $q$-systems is well defined up to the natural
notions of equivalence; see
\cite{Randrianarisoa2020,SheekeyVandeVoorde2020}.

We now recall the basic relation between rank weights and hyperplane
intersections.  We use the nondegenerate $\fqm$-bilinear form
$\sigma$ fixed above.  For $x\in\fqm^k\setminus\{0\}$, let
\[
x^\perp
=
\{u\in\fqm^k:\sigma(x,u)=0\}
\]
be the corresponding $\fqm$-hyperplane.  If the columns of $G$ form
an $\fq$-basis of $U$, then
\[
xG=
\bigl(\sigma(x,u_1),\ldots,\sigma(x,u_n)\bigr)
\]
up to the choice of coordinates associated with $\sigma$.  Therefore
the rank weight of $xG$ is the $\fq$-dimension of the image of the
map
\[
\varphi_x:U\longrightarrow\fqm,
\qquad
u\longmapsto\sigma(x,u).
\]
Since $\ker\varphi_x=U\cap x^\perp$,
the rank-nullity theorem gives
\begin{equation}\label{eq:code-system}
\wt(xG)
=
n-\dimq(U\cap x^\perp).
\end{equation}

Thus the rank-weight spectrum of $\cC$ is completely determined by
the hyperplane-intersection spectrum of the associated $q$-system
$U$.  More precisely, if
\[
N_j(U)
=
\#\bigl\{
H: H \,\,\text{hyperplane of}\,\,\PG(k-1,q^m) \,\,\text{and}\,\,
\dimq(U\cap H)=j
\bigr\},
\]
then, for $i>0$, the number $A_i(\cC)$ of codewords of rank weight
$i$ satisfies
\[
A_i(\cC)
=
(q^m-1)N_{n-i}(U).
\]
Indeed, the nonzero vectors defining the same projective point differ
by multiplication by an element of $\fqm^*$ and hence determine the
same hyperplane and codeword rank.

%Formula~\eqref{eq:code-system} remains valid more generally for any
%$\fq$-subspace $U\le\fqm^k$, with
%\[
%n=\dimq U,
%\]
%and for any matrix $G$ whose columns form an $\fq$-basis of $U$; in
%this more general situation, $U$ need not span the whole space over
%$\fqm$.

There is a second formulation of the correspondence which will be
particularly useful in this paper.  Combining
\eqref{eq:code-system} with
\eqref{eq:hyperplane-point-duality}, applied to
\[
H=x^\perp,
\qquad
H^\perp=\langle x\rangle_{\fqm},
\]
gives
\begin{equation}\label{eq:dual-point-code-weight}
w_{U^{\perp'}}(\langle x\rangle_{\fqm})
=
m-\wt(xG),
\end{equation}
see \cite{BorelloZullo2024}.
Hence the point weights of the dual subspace $U^{\perp'}$ are exactly
the \emph{deficiencies} of the rank weights of $\cC$ from $m$.
In particular,
\[
\wt(xG)<m
\quad\Longleftrightarrow\quad
\langle x\rangle_{\fqm}\in L_{U^{\perp'}}.
\]
Formulas \eqref{eq:code-system} and \eqref{eq:dual-point-code-weight} remain valid for an
arbitrary $\fq$-subspace $U\le\fqm^k$, with $n=\dimq U$, and for any matrix whose columns form an
$\fq$-basis of $U$; in this more general situation $U$ doesn't need to span $\fqm^k$ over $\fqm$.

The following consequence will be useful when studying few-weight
codes.  It determines exactly when the maximum possible rank weight
$m$ occurs.

\begin{proposition}
\label{prop:maximum-rank-weight}
Let $\cC$ be a nondegenerate $[n,k]_{q^m/q}$ rank-metric code, and
let $U$ be an associated $q$-system.  The following are equivalent:
\begin{enumerate}[label=\textup{(\alph*)}]
\item no codeword of $\cC$ has rank weight $m$;
\item $L_{U^{\perp'}}=\PG(k-1,q^m)$;
\item $n\le m-1$.
\end{enumerate}
Consequently, every nondegenerate $[n,k]_{q^m/q}$ code with
$n\ge m$ contains a codeword of rank weight $m$.
\end{proposition}

\begin{proof}
By \eqref{eq:dual-point-code-weight}, a point
$\langle x\rangle_{\fqm}$ has positive weight with respect to
$U^{\perp'}$ if and only if $\wt(xG)<m$.
Hence \textup{(a)} and \textup{(b)} are equivalent.

Put $N=\dimq U^{\perp'}=km-n$.
By Lemma~\ref{lem:covering}, applied with $S=\fqm^k$, one has
$L_{U^{\perp'}}=\PG(k-1,q^m)$
if and only if
$N>(k-1)m$.
Since $N=km-n$, this is equivalent to $n<m$,
or, equivalently, $n\le m-1$.
\end{proof}

The last assertion of Proposition~\ref{prop:maximum-rank-weight} is a special case of \cite[Proposition~3.11]{AlfaranoEtAl2022}, according to which every nondegenerate $[n,k]_{q^m/q}$ code contains a codeword of rank $\min\{m,n\}$.

\begin{remark}\label{rem:three-weight-maximum}
In particular, if a nondegenerate %three-weight 
rank-metric code has
length $n\ge m$, then one of its %three
nonzero weights is necessarily
$m$.  Under geometric duality, the weight $m$ corresponds to the
zero point weight of $U^{\perp'}$, whereas every nonmaximum codeword
weight $w<m$ corresponds to the positive point weight $m-w$ of $L_{U^{\perp'}}$.  Thus a code with nonzero weights
\[
d_1<\cdots<d_r<m
\qquad\text{and}\qquad m
\]
corresponds on the dual side to a linear set whose positive point
weights are
\[
m-d_r<\cdots<m-d_1.
\]
\end{remark}

We finally recall the operation of field expansion, which provides the
link between subfield linearity of a $q$-system and divisibility of
rank weights.

\begin{definition}\label{def:expansion}
Let $e\mid m$, and suppose that an $[n,k]_{q^m/q}$ system
$U\le\fqm^k$ is an $\mathbb F_{q^e}$-subspace.  Then necessarily
$e\mid n$, and
\[
\dim_{\mathbb F_{q^e}}U=\frac ne.
\]
Choose an $\mathbb F_{q^e}$-basis
$u_1,\ldots,u_{n/e}$ of $U$, and let $G_e$ be the matrix having these vectors as columns.
The code generated by $G_e$ is an
$[n/e,k]_{q^m/q^e}$ rank-metric code.

Fix an $\fq$-basis
$\beta_1,\ldots,\beta_e$ of $\mathbb F_{q^e}$.  Replacing each column
$u_i$ of $G_e$ by
$\beta_1u_i,\ldots,\beta_eu_i$
gives a matrix with $n$ columns whose $\fq$-span is $U$.  The
corresponding $[n,k]_{q^m/q}$ code is called the
$\fq$-\textbf{expansion} of the original code  $\cC_e$ generated by $G_e$  and is denoted by $\Exp_{q^e/q}(\cC_e)$.
\end{definition}

If $c$ is a codeword before expansion and $\widetilde c$ is the
corresponding codeword after expansion, then
\[
\wt_q(\widetilde c)
=
e\,\wt_{q^e}(c).
\]
Consequently, every rank weight of an $\fq$-expansion from
$\mathbb F_{q^e}$ is divisible by $e$.  On the geometric side, this
is reflected by the fact that, when $U$ is
$\mathbb F_{q^e}$-linear, every intersection
$U\cap H$
with an $\fqm$-subspace $H$ is again an
$\mathbb F_{q^e}$-subspace, and hence all the corresponding
intersection dimensions over $\fq$ are multiples of $e$.

\begin{remark}\label{rem:simplex-expansion}
A \textbf{simplex code} is a nondegenerate $[mk,k]_{q^m/q}$ code $\cC$.  Its associated system has $\fq$-dimension $mk$, hence it is the whole space $U=\fqm^k$, which is trivially $\fqm$-linear.  Therefore simplex codes are particular expanded codes: $\cC$ is, up to equivalence, the $\fq$-expansion $\Exp_{q^m/q}(\cC_m)$ of the $[k,k]_{q^m/q^m}$ code $\cC_m=\fqm^k$, whose generator matrix is the identity matrix.  Since every nonzero codeword of $\cC_m$ has rank weight one over $\fqm$, every nonzero codeword of $\cC$ has rank weight $m$; this agrees with \eqref{eq:code-system}, as $\dimq(\fqm^k\cap H)=(k-1)m$ for every hyperplane $H$.
\end{remark}

\section{Fields of linearity of few-weight linear sets}\label{sec:linear-sets}

We recall two results that are central to the paper.

\begin{theorem}\cite[Theorem~2]{CsajbokMarinoPepe2024}\label{thm:CMP}
Let $U$ define a linear set of rank $N$ spanning $\PG(h-1,q^m)$, and assume
\[
  N\le (h-1)m,\qquad m\le q.
\]
If every point of $L_U$ has weight at least $a\ge2$, then there exists a divisor $t$ of $m$, with $t\ge a$, such that
\[
  L_U=L_{\langle U\rangle_{\mathbb F_{q^t}}}.
\]
\end{theorem}

Repeated application produces a representation over a maximal subfield whose associated linear set contains a point of relative weight one.

\begin{theorem}\cite[Corollary~3.16]{AdriaensenSantonastaso2024}\label{thm:minsize}
Let $Q$ be a prime power and let $L_W$ be an $\mathbb F_Q$-linear set of rank $g\ge2$ spanning a projective subspace of positive dimension.  If $L_W$ contains a point of $\mathbb F_Q$-weight one, then
\[
  |L_W|\ge Q^{g-1}+1.
\]
\end{theorem}

The above result gives the stronger bound $|L_W|\ge Q^{g-1}+I$, where $I$ is the number of secant lines through a point of weight one.  Since the set spans at least a line, $I\ge1$.  See also \cite{DeBeuleVandeVoorde2019} for earlier results on the minimum size of linear sets.

The next known result treats ranks that are multiples of the extension degree.

\begin{proposition}\label{prop:multiple-m}
Let $U$ define a linear set whose rank is a multiple of $m$, and let $a$ be its minimum positive point weight.  Then $a\mid m$ and $U$ is $\mathbb F_{q^a}$-linear.  Consequently, $a$ divides every point weight and
\[
  a=\gcd\{w_U(P):P\in L_U\}.
\]
Moreover, $\mathbb F_{q^a}$ is the maximum field of linearity of $U$.
\end{proposition}

\begin{proof}
For $a>1$, the first assertion is \cite[Theorem~1]{CsajbokMarinoPepe2024}; for $a=1$ it is immediate.  Once $U$ is $\mathbb F_{q^a}$-linear, every intersection $U\cap P$ is an $\mathbb F_{q^a}$-subspace, so every point weight is divisible by $a$.  Since $a$ itself occurs, it equals the $\gcd\{w_U(P):P\in L_U\}$.  If $U$ is $\mathbb F_{q^s}$-linear, then $s$ divides every point weight, hence $s\mid a$ and $\mathbb F_{q^s}\subseteq\mathbb F_{q^a}$.
\end{proof}

We next isolate the elementary gcd consequence of the vector-counting identity.

\begin{lemma}\label{lem:gcd-rank}
Let $U$ be an $N$-dimensional $\fq$-subspace and assume that its positive point weights belong to $A=\{a_1,\ldots,a_r\}$.
Put $e=\gcd(a_1,\ldots,a_r)$.  Then $e\mid N$.
\end{lemma}

\begin{proof}
By \eqref{eq:vector-counting},
\[
  q^N-1=\sum_{i=1}^rN_{a_i}(q^{a_i}-1).
\]
Since $e\mid a_i$, the integer $q^e-1$ divides every summand on the right.  Hence $q^e-1\mid q^N-1$, which is equivalent to $e\mid N$.
\end{proof}

\begin{definition}\label{def:gap}
For positive integers $e<N$ define $\cD_m(e,N)=\{t:t\mid m,\ e<t\le N-e\}$.
We say that $(m,e,N)$ satisfies the \textbf{divisor-gap condition} if $\cD_m(e,N)=\varnothing$.
Equivalently, if
\[
 d^+(m,e)=\min\{t:t\mid m,\ t>e\},
\]
with $d^+(m,e)=+\infty$ when the set is empty, then the condition is
\[
  d^+(m,e)>N-e.
\]
\end{definition}

\begin{example}
Take $m=12$ and $e=3$.  If $N=6$, then
\[
\cD_{12}(3,6)
=
\{t:t\mid 12,\ 3<t\le3\}
=
\varnothing,
\]
so $(12,3,6)$ satisfies the divisor-gap condition.  Equivalently,
the smallest divisor of $12$ strictly larger than $3$ is
$d^+(12,3)=4$, and $4>N-e=3$.

On the other hand, if $N=9$, then
\[
\cD_{12}(3,9)
=
\{t:t\mid12,\ 3<t\le6\}
=
\{4,6\},
\]
so the divisor-gap condition does not hold.
\end{example}

We now combine the preceding results to show that a lower bound on the point weights not only forces a representation over a larger subfield, but also determines precisely the rank of such a representation.

\begin{theorem}\label{thm:larger-field-representation}
Let $q\ge m$, let $U\le_{\fq}V$, and put $S=\langle U\rangle_{\fqm}$, $h=\dim_{\fqm}S\ge 2$, $N=\dimq U$.
Assume that $N\le (h-1)m$,
and let
\[
    a=\min\{w_U(P):P\in L_U\}.
\]
If $a>1$, then there exists a divisor $t\mid m$, with $t\ge a$,
such that, setting
\[
    W=\langle U\rangle_{\mathbb F_{q^t}},
\]
one has $L_U=L_W$.
Moreover, $t$ may be chosen so that $L_W$, regarded as an $\mathbb F_{q^t}$-linear set, contains a point of $\mathbb F_{q^t}$-weight one.  For such a choice,
\[
    \dimq W=N+t-a.
\]
In particular,
\[
    t\mid N-a
    \qquad\text{and}\qquad
    a\le t\le N-a.
\]
\end{theorem}

\begin{proof}
Since $h\ge2$, $N\le(h-1)m$, $m\le q$ and every point of $L_U$ has weight at least $a\ge2$, Theorem~\ref{thm:CMP} applies.  Starting from $U$ and enlarging the field of linearity whenever the resulting representation has no point of relative weight one, we obtain a divisor $t\mid m$ and an $\mathbb F_{q^t}$-subspace
$W=\langle U\rangle_{\mathbb F_{q^t}}$
such that
\begin{equation}\label{eq:same-set}
  L_W=L_U,
\end{equation}
$t\ge a$, and $L_W$, regarded as an $\mathbb F_{q^t}$-linear set, has a point of $\mathbb F_{q^t}$-weight one.
Indeed, since $N\le(h-1)m$, Lemma~\ref{lem:covering} shows that $L_U\ne\PG(S,\fqm)$.  By the same lemma, every $\fq$-subspace of $S$ defining the point set $L_U$ has $\fq$-dimension at most $(h-1)m$; for an $\mathbb F_{q^t}$-subspace, this means $\mathbb F_{q^t}$-dimension at most $(h-1)m/t$, where $m/t$ is the degree of $\fqm$ over $\mathbb F_{q^t}$.  Furthermore,
\[
  \frac mt\le m\le q\le q^t,
\]
so Theorem~\ref{thm:CMP} can be reapplied over the base field $\mathbb F_{q^t}$ whenever the current representation has no point of relative weight one.  Each application strictly increases the degree of the field of linearity, which is a divisor of $m$, so the process terminates.

Write
\[
  g=\dim_{\mathbb F_{q^t}}W,
  \qquad N'=\dimq W=gt.
\]
Because $U\subseteq W$,
\begin{equation}\label{eq:Nprime-ge-N}
  N'\ge N.
\end{equation}
Moreover, $W$ spans $S$ over $\fqm$, and hence $g\ge h\ge2$.

Apply Theorem~\ref{thm:minsize} over the field $\mathbb F_{q^t}$.  We obtain
\begin{equation}\label{eq:lower-size}
  |L_U|=|L_W|\ge (q^t)^{g-1}+1=q^{N'-t}+1.
\end{equation}
On the other hand, every point of $L_U$ has weight at least $a$, so \eqref{eq:size-upper-minweight} yields
\begin{equation}\label{eq:upper-size}
  |L_U|\le \frac{q^N-1}{q^{a}-1}.
\end{equation}
Combining \eqref{eq:lower-size} and \eqref{eq:upper-size},
\[
  q^{N'-t}+1\le \frac{q^N-1}{q^{a}-1}<q^{N-a+1}.
\]
Since the exponents are integers, it follows that
\begin{equation}\label{eq:key-exponent}
  N'-t\le N-a.
\end{equation}

Now choose a point $P\in L_U$ of minimum weight, so that $w_U(P)=a$.
Since $L_U=L_W$, we have $P\in L_W$.  Moreover, $W$ is an $\mathbb F_{q^t}$-subspace, and therefore
\[
    w_W(P)=\rho t
\]
for some integer $\rho\ge1$.  Since $U\cap P\subseteq W\cap P$ and $(W\cap P)\cap U=U\cap P$, the quotient $(W\cap P)/(U\cap P)$ embeds into $W/U$.  Hence
\[
    N'-N
    \ge
    w_W(P)-w_U(P)
    =
    \rho t-a
    \ge
    t-a.
\]
Combining this with \eqref{eq:key-exponent}, which gives $N'-N\le t-a$, we obtain
\[
    N'=N+t-a.
\]
Finally, $N-a=N'-t=(g-1)t$.  Hence $t\mid N-a$, and $t\le N-a$ because $g\ge2$.
\end{proof}

\begin{remark}\label{rem:origin-dimension-formula}
The proof of Theorem~\ref{thm:larger-field-representation} extends the argument
of \cite[Corollary~III.3]{ZulloEtAl2024}, which concerns $\fq$-subspaces of
$\fqm^2$ all of whose points have weight $0$ or $e$. There, Theorem~\ref{thm:CMP}
is iterated until a point of relative weight one appears, and the bound of
\cite[Theorem~1.2]{DeBeuleVandeVoorde2019} is compared with the exact size
$(q^n-1)/(q^e-1)$ of $L_U$, giving $n-e=\dimq\langle U\rangle_{\mathbb F_{q^t}}-t$,
that is, the identity $\dimq W=N+t-a$ in that special case. Here only the minimum
weight $a$ is prescribed and $|L_U|$ is unknown: the exact count is replaced by
\eqref{eq:size-upper-minweight} together with the inequality $\dimq W-N\ge t-a$,
and Theorem~\ref{thm:minsize} allows us to work directly in $\PG(h-1,q^m)$
instead of restricting to lines as in \cite{ZulloEtAl2024}.
\end{remark}

The divisor-gap condition introduced above allows us to strengthen the larger-field representation of Theorem~\ref{thm:larger-field-representation} into a structural result, forcing the field of linearity to be determined by the gcd of the positive point weights.

\begin{theorem}\label{thm:gcd-rigidity}
Let $U\le_{\fq}V$ satisfy the assumptions of Theorem~\ref{thm:larger-field-representation}, and suppose that the positive point weights of $L_U$ are
\[
    A=\{a_1,\ldots,a_r\},
    \qquad
    1<a_1<\cdots<a_r.
\]
Put $e=\gcd(a_1,\ldots,a_r)$ and $N=\dimq U$.
Assume that $\cD_m(e,N)=\varnothing$.
Let $t\mid m$ and $W=\langle U\rangle_{\mathbb F_{q^t}}$ be as in Theorem~\ref{thm:larger-field-representation}, with $t$ chosen so that $L_W$ has a point of $\mathbb F_{q^t}$-weight one.  Then
\[
    a_1=t=e \,\,\,\text{and}\,\,\,    U=W.
\]
Consequently, $U$ is an $\mathbb F_{q^e}$-subspace of $V$.  In particular, $e\mid m$, $e\mid N$, the minimum positive point weight $a_1$ divides every $a_i$, and $\mathbb F_{q^e}$ is the maximum field of linearity of $U$.
\end{theorem}

\begin{proof}
Theorem~\ref{thm:larger-field-representation}, applied with $a=a_1$, gives
\begin{equation}\label{eq:t-interval}
  e\le a_1\le t\le N-a_1\le N-e.
\end{equation}
The integer $t$ divides $m$.  If $t>e$, then \eqref{eq:t-interval} gives $t\in\cD_m(e,N)$, contradicting the hypothesis.  Hence $t=e$, and then \eqref{eq:t-interval} forces
\[
  a_1=t=e.
\]
In particular, $e\mid m$.  Moreover,
\[
  \dimq W=N+t-a_1=N.
\]
Since $U\subseteq W$, we conclude that $U=W$.  Thus $U$ is $\mathbb F_{q^e}$-linear, and $e\mid N$ (this also follows from Lemma~\ref{lem:gcd-rank}).  Finally, if $U$ were linear over $\mathbb F_{q^s}$, then every intersection $U\cap P$ would be an $\mathbb F_{q^s}$-subspace, so $s$ would divide every point weight and hence $s\mid e$.  Therefore $\mathbb F_{q^e}$ is the maximum field of linearity.
\end{proof}

\begin{remark}\label{rem:refined-gap}
The proof of Theorem~\ref{thm:gcd-rigidity} uses the divisor-gap condition only to exclude the divisor $t$ produced by Theorem~\ref{thm:larger-field-representation}, and this divisor is more constrained than an arbitrary element of $\cD_m(e,N)$.
\begin{enumerate}[label=\textup{(\roman*)}]
  \item The divisor $t$ is a proper divisor of $m$.  Indeed, if $t=m$, then $W=\langle U\rangle_{\fqm}=S$, so $L_U=L_W=\PG(S,\fqm)$, contradicting Lemma~\ref{lem:covering} since $N\le(h-1)m$.
  \item The divisor $t$ satisfies $a_1\le t\le N-a_1$ and $t\mid N-a_1$.
\end{enumerate}
Consequently, the conclusions of Theorem~\ref{thm:gcd-rigidity} remain valid under the weaker hypothesis that no proper divisor $s$ of $m$ satisfies
\[
  a_1<s\le N-a_1
  \qquad\text{and}\qquad
  s\mid N-a_1.
\]
Indeed, under this hypothesis $t=a_1$, so $\dimq W=N$ and $U=W$ is $\mathbb F_{q^{a_1}}$-linear; hence $a_1$ divides every point weight, and $a_1=e$.  By~(i), the divisor $m$ itself never obstructs; this is relevant only when $N\ge m+e$, which, as $N\le(h-1)m$, requires $h\ge3$.  These refinements are not used elsewhere in the paper.
\end{remark}

We now record some consequences.  

\begin{corollary}\label{cor:three-point}
Let $U$ satisfy the assumptions of Theorem~\ref{thm:larger-field-representation}, and suppose that the point weights are
\[
  0,s_1,s_2,
  \qquad 1<s_1<s_2.
\]
Put $e=\gcd(s_1,s_2)$ and $N=\dimq U$.  If $\cD_m(e,N)=\varnothing$, then
\[
  e=s_1,
\]
so
\[
  s_1\mid m,\qquad s_1\mid N,\qquad s_1\mid s_2,
\]
and $U$ is $\mathbb F_{q^{s_1}}$-linear.
\end{corollary}

\begin{corollary}\label{cor:prime-m}
Let $m$ be prime and $q\ge m$.  Let $U\le_{\fq}V$ be such that $S=\langle U\rangle_{\fqm}$ has dimension $h\ge2$, and assume that $\dimq U\le(h-1)m$.  Then $L_U$ has a point of weight one. 
\end{corollary}

\begin{proof}
Assume that every point of $L_U$ has weight at least two.  By Theorem~\ref{thm:CMP}, there is a divisor $t\ge2$ of $m$ with $L_U=L_{\langle U\rangle_{\mathbb F_{q^t}}}$.  Since $m$ is prime, $t=m$, and then $\langle U\rangle_{\fqm}=S$ gives $L_U=\PG(S,\fqm)$.  This contradicts Lemma~\ref{lem:covering}.
\end{proof}

\begin{remark}
The proof of Corollary~\ref{cor:prime-m} uses only Theorem~\ref{thm:CMP} and Lemma~\ref{lem:covering}, and it holds in arbitrary projective dimension.
\end{remark}

The preceding results admit a natural dual formulation in terms of hyperplane weights. Indeed, trace duality identifies the excess of a hyperplane weight above the minimum possible value \(n-m\) with the corresponding point weight of the dual linear set. This allows us to translate the above structural result directly into a criterion for the hyperplane-weight spectrum.

\begin{theorem}\label{thm:hyperplane-criterion}
Let $q\ge m$, let $U\le_{\fq}V$ have rank $n$, and put
\[
  W=U^{\perp'},\qquad N=\dimq W=km-n.
\]
Assume that $W$ spans an $\fqm$-subspace $S$ of dimension $h\ge2$ and that $N\le(h-1)m$.  Assume that the hyperplane weights of $U$ are precisely
\[
  n-m,\ n-m+a_1,\ldots,n-m+a_r,
  \qquad 1\le a_1<\cdots<a_r,
\]
and put $e=\gcd(a_1,\ldots,a_r)$.
If $\cD_m(e,N)=\varnothing$,
then $a_1=e$,
and both $W$ and $U$ are $\mathbb F_{q^e}$-linear.  In particular,
\[
  e\mid m,\qquad e\mid n,
\]
and every hyperplane weight of $U$ is divisible by $e$.
\end{theorem}

\begin{proof}
By Proposition~\ref{prop:weight-duality}, the point weights of $L_W$ are
\[
  0,a_1,\ldots,a_r.
\]
If $a_1=1$, then $e=1$ and there is nothing to prove.  Otherwise, Theorem~\ref{thm:gcd-rigidity} applies to $W$, gives $a_1=e$, and proves that $W$ is $\mathbb F_{q^e}$-linear.  By Lemma~\ref{lem:perp-linear}, $U=W^{\perp'}$ is also $\mathbb F_{q^e}$-linear.  Hence $e\mid\dimq U=n$, and every intersection $U\cap H$ with an $\fqm$-hyperplane $H$ is an $\mathbb F_{q^e}$-subspace.
\end{proof}

\begin{remark}\label{rem:dual-span}
The hypothesis that $U^{\perp'}$ spans a projective subspace of positive dimension is automatic as soon as at least two distinct positive excesses $a_i$ occur.  Indeed, if $U^{\perp'}$ were contained in one $\fqm$-point, then its linear set would have only one positive point weight, namely $N$.  Equivalently, $U^{\perp'}$ is contained in a point precisely when $U$ contains an $\fqm$-hyperplane.  This exceptional case should be separated before applying Theorem~\ref{thm:hyperplane-criterion}.

If $(k-1)m<n<km$ and $h\ge2$, then $N<m\le(h-1)m$, so the rank hypothesis of Theorem~\ref{thm:hyperplane-criterion} is automatic.
\end{remark}

\begin{corollary}\label{cor:three-hyperplanes}
Let the assumptions of Theorem~\ref{thm:hyperplane-criterion} hold and suppose that the hyperplane weights are
\[
  n-m,\qquad n-m+s_1,\qquad n-m+s_2,
  \qquad 1<s_1<s_2.
\]
Set $e=\gcd(s_1,s_2)$.  If $\cD_m(e,km-n)=\varnothing$,
then $e=s_1$, $U$ is $\mathbb F_{q^{s_1}}$-linear, and
\[
  s_1\mid m,\qquad s_1\mid n,\qquad s_1\mid s_2.
\]
\end{corollary}

\section{Antipodal $2$-dimensional codes}\label{sec:antipodal}

We call a rank-metric code of length $n$ \textbf{antipodal} if it contains a codeword of rank weight $n$.  If $\cC$ is a nondegenerate $[n,2]_{q^m/q}$ code with $0<n\le m$ and $U$ is an associated system, then Lemma~\ref{lem:covering}, with $h=2$, gives $L_U\ne\PG(1,q^m)$; if $x^\perp$ is a point outside $L_U$, then \eqref{eq:code-system} shows that $xG$ has rank $n$.  Hence every nondegenerate $[n,2]_{q^m/q}$ code with $n\le m$ is antipodal.
More generally, by \cite[Proposition~3.11]{AlfaranoEtAl2022}, every nondegenerate rank-metric code of length $n\le m$ is antipodal, whatever its dimension; thus, below, the word \emph{antipodal} is only a reminder of this property and imposes no extra condition.

Throughout this section, $\cC$ is a nondegenerate $[n,2]_{q^m/q}$ code with $n\le m$, $G$ is a generator matrix of $\cC$ and $U\le\fqm^2$ is the associated system.  We write the nonzero rank weights of $\cC$ as
\[
  d_1<\cdots<d_r<n.
\]
Note that $r\ge1$, since $n>0$. For $k=2$ the hyperplanes of $\PG(1,q^m)$ are points, and \eqref{eq:code-system} reads
\begin{equation}\label{eq:antipodal-point-weight}
  w_U(x^\perp)=n-\wt(xG),
  \qquad x\in\fqm^2\setminus\{0\}.
\end{equation}
Hence the positive point weights of $L_U$ are $n-d_r<\cdots<n-d_1$, and the minimum positive point weight of $L_U$ is
\[
  a=n-d_r,
\]
the gap between the maximum weight and the second largest weight of $\cC$.

For $i\ge1$ we denote by $A_i(\cC)$ the number of codewords of $\cC$ of rank weight $i$, and by $N_j(U)$ the number of points of $\PG(1,q^m)$ of weight $j$ with respect to $U$. Since the hyperplanes of $\PG(1,q^m)$ are its points, this agrees with the notation $N_j(U)$ of Section~\ref{sec:preliminaries}.  By \eqref{eq:antipodal-point-weight}, and since proportional nonzero vectors $x$ give the same point $x^\perp$, we have $A_i(\cC)=(q^m-1)N_{n-i}(U)$ for $1\le i\le n$.  Counting the nonzero codewords and applying \eqref{eq:vector-counting} to $U$, we get, for every nondegenerate $[n,2]_{q^m/q}$ code,
\begin{equation}\label{eq:antipodal-counting}
  \sum_{i=1}^{n}A_i(\cC)=q^{2m}-1,
  \qquad
  \sum_{i=1}^{n}A_i(\cC)\bigl(q^{\,n-i}-1\bigr)=(q^m-1)(q^n-1).
\end{equation}
These two relations will be used to compare the weight distribution of $\cC$ with the weight distributions of the codes obtained from it by puncturing or expansion.

Recall the $\fq$-expansion of Definition~\ref{def:expansion}.  Let $t\mid m$, let $\cC_t$ be an $\bigl[\tfrac nt,2\bigr]_{q^m/q^t}$ code and let $\cC=\Exp_{q^t/q}(\cC_t)$.  If $\delta_1<\cdots<\delta_r$ are the nonzero weights of $\cC_t$ smaller than its length $\tfrac nt$, then, since expansion multiplies every rank weight by $t$, the nonzero weights of $\cC$ are
\[
  d_i=t\,\delta_i\quad(i=1,\ldots,r)
  \qquad\text{and}\qquad
  n .
\]
In particular,
\[
  n-d_r=t\Bigl(\frac nt-\delta_r\Bigr)\ge t,
\]
so that, if $t>1$ (that is, if the expansion is not trivial), the gap between the maximum weight and the second largest weight of $\cC$ is at least $2$.  The next result shows that, for $q\ge m$, this necessary condition is sufficient up to puncturing, i.e.  every antipodal $2$-dimensional code whose gap is at least $2$ is obtained by deleting few coordinates from an $\fq$-expansion.

\begin{theorem}\label{thm:antipodal-puncturing}
Let $q\ge m$ and let $\cC$ be a nondegenerate $[n,2]_{q^m/q}$ code, with $0<n\le m$, whose nonzero weights are $d_1<\cdots<d_r<n$.  Assume that
\[
  a=n-d_r\ge2 .
\]
Then there exists a divisor $t$ of $m$ with
\[
  a\le t\le\min\{d_r,\,m-d_r\}
  \qquad\text{and}\qquad
  t\mid d_r,
\]
and a nondegenerate antipodal $[d_r+t,2]_{q^m/q}$ code $\cC'$, where $d_r+t=n+t-a\le m$, with the following properties.
\begin{enumerate}[label=\textup{(\roman*)}]
  \item $\cC'$ is, up to equivalence, the $\fq$-expansion $\Exp_{q^t/q}(\cC_t)$ of a $\bigl[\tfrac{d_r}{t}+1,2\bigr]_{q^m/q^t}$ code $\cC_t$.
  \item Up to equivalence, $\cC$ is obtained from $\cC'$ by deleting $t-a$ coordinates: the map
  \[
    \pi:\cC'\longrightarrow\cC
  \]
  deleting the last $t-a$ coordinates is an $\fqm$-linear bijection.
  \item For every $c'\in\cC'$,
  \[
    \frac{\wt(c')}{t}
    =1+\left\lfloor\frac{\wt(\pi(c'))-a}{t}\right\rfloor .
  \]
  \item For every $c'\in\cC'$,
  \[
    \wt(c')-(t-a)\le\wt(\pi(c'))\le\wt(c').
  \]
  Consequently $\wt(c')=t\lceil\wt(\pi(c'))/t\rceil$, every nonzero weight $w$ of $\cC$ satisfies
  \[
    (w\bmod t) \in\{0\}\cup\{a,a+1,\ldots,t-1\},
  \]
  and $\wt(\pi(c'))=n$ if and only if $\wt(c')=d_r+t$.
  \item { For $i\ge1$, let $A_i(\cC)$, $A_i(\cC')$ and $A_i(\cC_t)$ be the numbers of codewords of rank weight $i$ in $\cC$, $\cC'$ and $\cC_t$, respectively \textup(the weight of $\cC_t$ being computed over $\mathbb F_{q^t}$\textup), and put $A_i(\cC)=0$ for $i>n$.  Then $A_i(\cC')=0$ whenever $t\nmid i$ and, for every $j=1,\ldots,\ell$, where $\ell=\frac{d_r}{t}+1$,
  \begin{equation}\label{eq:(v)1}
    A_{jt}(\cC')=A_j(\cC_t)=\sum_{i=(j-1)t+a}^{jt}A_i(\cC).
  \end{equation}
  In particular, $A_{d_r+t}(\cC')=A_n(\cC)$, $A_i(\cC)\le A_{t\lceil i/t\rceil}(\cC')$ for every $i\ge1$, and the weight distribution of $\cC$ satisfies, besides \eqref{eq:antipodal-counting}, the relation
  \begin{equation}\label{eq:(v)2}
\sum_{i=1}^{n}A_i(\cC)\bigl(q^{\,d_r+t-t\lceil i/t\rceil}-1\bigr)=(q^m-1)\bigl(q^{d_r+t}-1\bigr).
  \end{equation}}
  \item If $\cC'$ has $s$ distinct nonzero weights, then $\cC$ has at most $(s-1)(t-a+1)+1$ distinct nonzero weights.
\end{enumerate}
\end{theorem}

\begin{proof}
\emph{Step 1: a larger field of linearity.}
Let $U$ be an $\F_q$-system associated with $\cC$, so $\langle U\rangle_{\fqm}=\fqm^2$ and $ \dim_{\fq} U=n\leq m$.  By \eqref{eq:antipodal-point-weight}, every point of $L_U$ has weight at least $a\ge2$, and $a$ is attained.  Since $q\geq m$ and $n\leq m$, Theorem~\ref{thm:larger-field-representation} applies: there exists a divisor $t\mid m$ such that
\[
  W=\langle U\rangle_{\mathbb F_{q^t}}
\]
satisfies $L_W=L_U$ and
\[
  \dimq W=n+t-a=d_r+t,
  \qquad
  a\le t\le n-a=d_r,
  \qquad
  t\mid n-a=d_r .
\]
Since $n\le m$, Lemma~\ref{lem:covering} gives $L_W=L_U\ne\PG(1,q^m)$, and applying the same lemma to $W$ we get $\dimq W\le m$, that is $t\le m-d_r$.  Put
\[
  \ell=\dim_{\mathbb F_{q^t}}W=\frac{d_r+t}{t}=\frac{d_r}{t}+1,
  \qquad\text{so that}\qquad
  n=(\ell-1)t+a .
\]

\emph{Step 2: the weights of $W$.}
We claim that, for every point $Q$ of $\PG(1,q^m)$,
\begin{equation}\label{eq:W-ceiling}
  w_W(Q)=t\left\lceil\frac{w_U(Q)}{t}\right\rceil ,
\end{equation}
where $w_W(Q)=\dim_q (W\cap Q)$.
If $Q\notin L_U=L_W$, both sides vanish.  Otherwise, $U\cap Q\subseteq W\cap Q$, so $(W\cap Q)/(U\cap Q)$ embeds into $W/U$, which has $\fq$-dimension $t-a$.  Therefore
\[
  w_U(Q)\le w_W(Q)\le w_U(Q)+t-a .
\]
Since $W$ is $\mathbb F_{q^t}$-linear and $Q$ is an $\fqm$-subspace, $W\cap Q$ is an $\mathbb F_{q^t}$-subspace, and $w_W(Q)$ is a multiple of $t$.  The least multiple of $t$ which is at least $w_U(Q)$ lies between $w_U(Q)$ and $w_W(Q)$; as the interval $[w_U(Q),w_U(Q)+t-a]$ has length $t-a<t$, it contains at most one multiple of $t$, and \eqref{eq:W-ceiling} follows.

\emph{Step 3: construction of $\cC'$ and proof of \textup{(i)} and \textup{(ii)}.}
Let $u_1,\ldots,u_n$ be an $\fq$-basis of $U$ and complete it to an $\fq$-basis
\[
  u_1,\ldots,u_n,w_{n+1},\ldots,w_{n+t-a}
\]
of $W$.  Let $G'\in\fqm^{2\times(n+t-a)}$ be the matrix having these vectors as columns, and let $\cC'$ be the code generated by $G'$.  Its columns are $\fq$-linearly independent and span $\fqm^2$ over $\fqm$, so $\cC'$ is a nondegenerate $[n+t-a,2]_{q^m/q}$ code with associated system $W$; since $n+t-a\le m$, it is antipodal.  Since  $u_1,\ldots,u_n$ is an $\fq$-basis  of $U$, up to equivalence, the matrix $G$ consisting of the first $n$ columns of $G'$  is a generator matrix of $\cC$.  Then the map $\pi$ deleting the last $t-a$ coordinates satisfies
\[
  \pi(xG')=xG,\qquad x\in\fqm^2,
\]
and it is an $\fqm$-linear bijection, since both $x\mapsto xG'$ and $x\mapsto xG$ are isomorphisms from $\fqm^2$ onto $\cC'$ and $\cC$, respectively.  This proves~(ii).  Moreover, choosing an $\mathbb F_{q^t}$-basis of $W$ and expanding it with respect to a fixed $\fq$-basis of $\mathbb F_{q^t}$ as in Definition~\ref{def:expansion}, we obtain another $\fq$-basis of $W$; thus $\cC'$ is equivalent to $\Exp_{q^t/q}(\cC_t)$, where $\cC_t$ is the $[\ell,2]_{q^m/q^t}$ code generated by the matrix whose columns are the chosen $\mathbb F_{q^t}$-basis of $W$.  This proves~(i).

\emph{Step 4: proof of \textup{(iii)}.}
Let $x\in\fqm^2\setminus\{0\}$, $c'=xG'$ and $c=\pi(c')=xG$.  Applying \eqref{eq:code-system} to $U$ and to $W$, and using \eqref{eq:W-ceiling} with $Q=x^\perp$, we get
\[
  \wt(c)=n-w_U(x^\perp),
  \qquad
  \wt(c')=(n+t-a)-w_W(x^\perp)
  =\ell t-t\left\lceil\frac{n-\wt(c)}{t}\right\rceil .
\]
Since $n=(\ell-1)t+a$,
\[
  \left\lceil\frac{n-\wt(c)}{t}\right\rceil
  =(\ell-1)+\left\lceil\frac{a-\wt(c)}{t}\right\rceil
  =(\ell-1)-\left\lfloor\frac{\wt(c)-a}{t}\right\rfloor ,
\]
and therefore
\[
  \wt(c')=t+t\left\lfloor\frac{\wt(c)-a}{t}\right\rfloor ,
\]
which is~(iii).

\emph{Step 5: proof of \textup{(iv)}.}
Deleting coordinates of a vector cannot increase its rank weight, hence $\wt(\pi(c'))\le\wt(c')$.  On the other hand, by~(iii),
\[
  \frac{\wt(c')}{t}\le1+\frac{\wt(\pi(c'))-a}{t},
\]
that is, $\wt(c')-(t-a)\le\wt(\pi(c'))$.  Now $\wt(c')=jt$ for some integer $j\ge0$, by~(i).  Then $\wt(\pi(c'))\in[(j-1)t+a,\,jt]$, and since $1\le a\le t$ this interval contains no multiple of $t$ other than $jt$ and is contained in $\bigl((j-1)t,jt\bigr]$; hence $\lceil\wt(\pi(c'))/t\rceil=j$, and $(\wt(\pi(c'))\bmod t)\in\{0\}\cup\{a,\ldots,t-1\}$.  Since $\pi$ is a bijection, every nonzero weight of $\cC$ arises in this way.  Finally, by~(iii), $\wt(c')=\ell t=d_r+t$ if and only if $\lfloor(\wt(\pi(c'))-a)/t\rfloor=\ell-1=d_r/t$, that is, if and only if $\wt(\pi(c'))\ge d_r+a=n$, i.e.\ $\wt(\pi(c'))=n$.

{
\emph{Step 6: proof of \textup{(v)}.}
By~(i), every weight of $\cC'$ is a multiple of $t$, so $A_i(\cC')=0$ if $t\nmid i$.  Since equivalent codes have the same weight distribution and the expansion map $\cC_t\to\Exp_{q^t/q}(\cC_t)$ is a bijection, multiplying every rank weight by $t$, we have $A_{jt}(\cC')=A_j(\cC_t)$.  By~(iv), $\wt(c')=t\lceil\wt(\pi(c'))/t\rceil$ for every $c'\in\cC'$.  Hence, for $1\le j\le\ell$, the map $\pi$ sends the codewords of $\cC'$ of weight $jt$ bijectively onto the codewords $c\in\cC$ with $(j-1)t<\wt(c)\le jt$, and by the residue condition in~(iv) the weights of these codewords lie in $[(j-1)t+a,\,jt]$.  This gives \eqref{eq:(v)1}.  For $j=\ell$ the interval is $[n,\,d_r+t]$, because $(\ell-1)t+a=d_r+a=n$, and therefore $A_{d_r+t}(\cC')=A_n(\cC)$.  The inequality $A_i(\cC)\le A_{t\lceil i/t\rceil}(\cC')$ is clear, since $A_i(\cC)$ is one of the summands of the formula with $j=\lceil i/t\rceil$.  Finally, $\cC'$ is a nondegenerate $[d_r+t,2]_{q^m/q}$ code, so the second relation in \eqref{eq:antipodal-counting} applied to $\cC'$ gives
\[
  \sum_{j=1}^{\ell}A_{jt}(\cC')\bigl(q^{\,d_r+t-jt}-1\bigr)=(q^m-1)\bigl(q^{d_r+t}-1\bigr).
\]
Replacing $A_{jt}(\cC')$ by $\sum_{\lceil i/t\rceil=j}A_i(\cC)$ we obtain \eqref{eq:(v)2}.

\emph{Step 7: proof of \textup{(vi)}.}
By~(v), the codewords of $\cC'$ of weight $jt$ are mapped by $\pi$ to codewords of $\cC$ whose weights lie in $\{jt-(t-a),\ldots,jt\}$, a set of $t-a+1$ integers, while the codewords of maximum weight $d_r+t$ are mapped exactly onto the codewords of weight $n$.  Since $\pi$ is a bijection, $\cC$ has at most $(s-1)(t-a+1)+1$ distinct nonzero weights.
}
\end{proof}

\begin{remark}\label{rem:antipodal-consequences}
We collect some consequences of Theorem~\ref{thm:antipodal-puncturing}.
\begin{enumerate}[label=\textup{(\alph*)}]
  \item Every nonzero codeword $c'\in\cC'$ has weight at least $t$, so by Theorem~\ref{thm:antipodal-puncturing}(iv) every nonzero codeword of $\cC$ has weight at least $t-(t-a)=a$.  Thus $d_1\ge n-d_r$, i.e.\ $d_1+d_r\ge n$.
  \item If $t=a$, then $\dimq W=n$, so $W=U$, no coordinate is deleted and $\cC$ itself is, up to equivalence, $\Exp_{q^a/q}(\cC_a)$ for an $\bigl[\tfrac na,2\bigr]_{q^m/q^a}$ code $\cC_a$.  In particular $a\mid m$, $a\mid n$ and $a$ divides every weight of $\cC$, so that $a=\gcd(n-d_1,\ldots,n-d_r)$.
  \item If $n=m$, then $t\le m-d_r=a$, and hence $t=a$.  So, in the square case, $\cC$ is an $\fq$-expansion from $\mathbb F_{q^a}$, in accordance with Remark~\ref{rem:square-case}.
  \item More generally, $t=a$ whenever no divisor $s$ of $m$ satisfies
  \[
    a<s\le\min\{d_r,\,m-d_r\}
    \qquad\text{and}\qquad
    s\mid d_r .
  \]
  If $e=\gcd(n-d_1,\ldots,n-d_r)$, then $e\le a$ and $d_r=n-a\le n-e$; hence this condition is implied by $\cD_m(e,n)=\varnothing$.  Therefore, if $\cD_m(e,n)=\varnothing$, then $t=a=e$ and, by~(b), $\cC$ is, up to equivalence, the $\fq$-expansion of an $[n/e,2]_{q^m/q^e}$ code.  For $n=m$ the same conclusion holds without any assumption on $q$, by Theorem~\ref{thm:length-multiple-m}; see Remark~\ref{rem:square-case}.
\end{enumerate}
\end{remark}

We now specialize to three-weight codes.

\begin{corollary}\label{cor:antipodal-puncturing-three}
Let $q\ge m$ and let $\cC$ be a nondegenerate $[n,2]_{q^m/q}$ code, with $n\le m$, whose nonzero weights are $d<w<n$, and assume that $n-w\ge2$.  Then there exists a divisor $t$ of $m$ with
\[
  n-w\le t\le\min\{w,\,m-w\}
  \qquad\text{and}\qquad
  t\mid w,
\]
such that the following hold.
\begin{enumerate}[label=\textup{(\roman*)}]
  \item $d\ge n-w$.
  \item Either $t\mid d$, or $(d\bmod t)\ge n-w$.
  \item $\cC$ is obtained, up to equivalence, by deleting $t-n+w$ coordinates from $\Exp_{q^t/q}(\cC_t)$, where $\cC_t$ is a $\bigl[\tfrac wt+1,2\bigr]_{q^m/q^t}$ code whose nonzero weights are the distinct values among
  \[
    \left\lceil\frac dt\right\rceil,\qquad \frac wt,\qquad \frac wt+1 .
  \]
  More precisely, the puncturing map sends the codewords of $\Exp_{q^t/q}(\cC_t)$ of weight $w+t$ bijectively onto the codewords of $\cC$ of weight $n$; { if $d\le w-t$, it sends those of weight $t\lceil d/t\rceil$ and $w$ bijectively onto those of weight $d$ and $w$, respectively, whereas if $d>w-t$ it sends the codewords of weight $w$ bijectively onto the codewords of weight $d$ or $w$}.  In particular, $\cC_t$ has at most three nonzero weights, and exactly two if $d>w-t$.
  \item Let $\cC'=\Exp_{q^t/q}(\cC_t)$ and put
  \[
    E=w+t-t\left\lceil\frac dt\right\rceil .
  \]
  The weight distributions of $\cC$, $\cC'$ and $\cC_t$ are related by $A_{jt}(\cC')=A_j(\cC_t)$ for every $j$, $A_{w+t}(\cC')=A_n(\cC)$, and
  \[
    \begin{cases}
      A_{t\lceil d/t\rceil}(\cC')=A_d(\cC),\quad A_w(\cC')=A_w(\cC), & \text{if } d\le w-t,\\[3pt]
      A_w(\cC')=A_d(\cC)+A_w(\cC), & \text{if } d>w-t .
    \end{cases}
  \]
  If moreover $t>n-w$, that is, if at least one coordinate is deleted, then
  \[
    \Delta=(q^{n-d}-1)(q^t-1)-(q^{n-w}-1)(q^E-1)\ne0
  \]
  and the weight distribution of $\cC$ is determined by $q,m,n,d,w$ and $t$:
  \begin{align*}
    A_d(\cC)&=\frac{q^m-1}{\Delta}\Bigl[(q^n-1)(q^t-1)-(q^{n-w}-1)(q^{w+t}-1)\Bigr],\\
    A_w(\cC)&=\frac{q^m-1}{\Delta}\Bigl[(q^{n-d}-1)(q^{w+t}-1)-(q^E-1)(q^n-1)\Bigr],\\
    A_n(\cC)&=q^{2m}-1-A_d(\cC)-A_w(\cC).
  \end{align*}
  In particular, if $d>w-t$, then $E=t$ and $A_d(\cC)+A_w(\cC)=(q^m-1)\dfrac{q^{w+t}-1}{q^t-1}$.
  \item If $n=m$, or more generally if no divisor $s$ of $m$ satisfies $n-w<s\le\min\{w,m-w\}$ and $s\mid w$, then $t=n-w$ and
  \[
    n-w\mid m,\qquad n-w\mid n,\qquad n-w\mid d,\qquad n-w\mid w;
  \]
  moreover, $\cC$ is, up to equivalence, the $\fq$-expansion of an $\bigl[\tfrac n{n-w},2\bigr]_{q^m/q^{n-w}}$ code with nonzero weights
  \[
    \frac{d}{n-w},\qquad\frac{w}{n-w},\qquad\frac{n}{n-w}.
  \]
\end{enumerate}
\end{corollary}

\begin{proof}
Apply Theorem~\ref{thm:antipodal-puncturing} with $r=2$, $d_1=d$, $d_2=w$ and $a=n-w$.  Assertion~(i) is Remark~\ref{rem:antipodal-consequences}(a), and~(ii) is the residue condition of Theorem~\ref{thm:antipodal-puncturing}(iv) applied to $d$.  For~(iii), Theorem~\ref{thm:antipodal-puncturing}(iv) shows that a codeword $c'$ with $\pi(c')$ of weight $d$, $w$ or $n$ has weight $t\lceil d/t\rceil$, $t\lceil w/t\rceil=w$ or $t\lceil n/t\rceil=w+t$, respectively, the last equality because $n=w+a$ with $t\mid w$ and $1\le a\le t$.  Since $\pi$ is a bijection and expansion multiplies weights by $t$, the weights of $\cC_t$ are as stated.  As $\lceil d/t\rceil\le w/t$, with equality if and only if $d>w-t$, the last claim of~(iii) follows.

For~(iv), recall that the nonzero weights of $\cC$ are $d$, $w$ and $n=w+a$, with $a=n-w$.  By Theorem~\ref{thm:antipodal-puncturing}(v), $A_{jt}(\cC')=A_j(\cC_t)$ is the number of codewords of $\cC$ whose weight lies in $[(j-1)t+a,jt]$.  The weight $d$ lies in the block $j=\lceil d/t\rceil$, the weight $w$ in the block $j=w/t$ and the weight $n$ in the block $j=w/t+1$; the first two blocks coincide if and only if $\lceil d/t\rceil=w/t$, that is, if and only if $d>w-t$.  This proves the relations between the weight distributions.  Now put $\alpha=A_d(\cC)/(q^m-1)$, $\beta=A_w(\cC)/(q^m-1)$ and $D=n-d$.  The second relation in \eqref{eq:antipodal-counting} and the last relation in Theorem~\ref{thm:antipodal-puncturing}(v), where $d_r=w$, $t\lceil w/t\rceil=w$ and $t\lceil n/t\rceil=w+t$, give the linear system
\[
  \alpha(q^D-1)+\beta(q^a-1)=q^n-1,
  \qquad
  \alpha(q^E-1)+\beta(q^t-1)=q^{w+t}-1,
\]
whose determinant is $\Delta$.  Assume $t>a$.  Since $d<w$ and $\lceil d/t\rceil\le w/t$, we have $D>a$ and $E\ge t>a$.  If $\Delta=0$, then expanding and dividing by $q^a$ we would get
\[
  q^{D+t-a}-q^{D-a}-q^{t-a}=q^E-q^{E-a}-1,
\]
where all the exponents are positive; the left-hand side is divisible by $q$, while the right-hand side is congruent to $-1$ modulo $q$, a contradiction.  Hence $\Delta\ne0$, and Cramer's rule gives $A_d(\cC)$ and $A_w(\cC)$; $A_n(\cC)$ follows from the first relation in \eqref{eq:antipodal-counting}.  If $d>w-t$, then $t\lceil d/t\rceil=w$, so $E=t$ and the second equation reads $(\alpha+\beta)(q^t-1)=q^{w+t}-1$.  Note that, when $t=a$, the two equations coincide, since $E=D$ and $\cC=\cC'$.

Finally, under the hypotheses of~(v), Remark~\ref{rem:antipodal-consequences}(c),(d) gives $t=n-w$; then~(ii) forces $t\mid d$, because $\{a,\ldots,t-1\}=\varnothing$, and Remark~\ref{rem:antipodal-consequences}(b) gives the remaining assertions.
\end{proof}

\begingroup
\subsection{Families of three-weight antipodal $2$-dimensional codes}\label{subsec:antipodal-examples}

In this subsection $\cC$ is a nondegenerate three-weight $[n,2]_{q^m/q}$ code, with $n\le m$ and nonzero weights $d<w<n$, and $U$ is its associated $q$-system.  We call
\[
  g_1=n-w
  \qquad\text{and}\qquad
  g_2=w-d
\]
the \emph{first gap} (between the maximum and the second largest weight) and the \emph{second gap} (between the second largest and the minimum weight) of $\cC$.  By Corollary~\ref{cor:antipodal-puncturing-three}, if $q\ge m$ and $g_1\ge2$, the first gap forces strong arithmetic restrictions: there is a divisor $t$ of $m$ with $g_1\le t\le\min\{w,m-w\}$ and $t\mid w$, so that in particular $\gcd(w,m)\ge2$, and if $t>g_1$ the weight distribution of $\cC$ is determined by its parameters.  We now describe some families of examples, which exist for every prime power $q$, showing different behaviours of the two gaps.  By \eqref{eq:antipodal-point-weight}, constructing $\cC$ amounts to constructing an $n$-dimensional $\fq$-subspace $U$ of $\fqm^2$, with $\langle U\rangle_{\fqm}=\fqm^2$, whose point weights take exactly the three values $0<n-w<n-d$; the weight distribution of $\cC$ is then given by $A_i(\cC)=(q^m-1)N_{n-i}(U)$.

We shall use the following terminology.  Let $t\mid m$.  An $\mathbb F_{q^t}$-subspace $U'$ of $\fqm^2$ is \textbf{$\mathbb F_{q^t}$-scattered} if $\dim_{\mathbb F_{q^t}}(U'\cap P)\le1$ for every point $P$ of $\PG(1,q^m)$; for $t=1$ we simply say that $U'$ is \textbf{scattered}, that is, all its points have weight at most one.  For $i\ge2$, an \textbf{$i$-club} of rank $N$ is an $\fq$-linear set $L_{U'}\subseteq\PG(1,q^m)$ of rank $N$ having exactly one point of weight $i$, called its \emph{head}, all the other points of $L_{U'}$ having weight one.

For every $t\mid m$, the $\mathbb F_{q^t}$-subspace
\[
  \{(x,x^{q^t}):x\in\fqm\}
\]
of $\mathbb F_{q^m}^2$ is $\mathbb F_{q^t}$-scattered: if $(x,x^{q^t})$ and $(y,y^{q^t})$, with $xy\ne0$, define the same point, then $(y/x)^{q^t-1}=1$, i.e.\ $y/x\in\mathbb F_{q^t}$.  Since an $\mathbb F_{q^t}$-subspace of an $\mathbb F_{q^t}$-scattered subspace is again $\mathbb F_{q^t}$-scattered, for every $\ell\le m/t$ there exist $\mathbb F_{q^t}$-scattered subspaces of $\fqm^2$ of $\mathbb F_{q^t}$-dimension $\ell$.
\begin{example}[Family $\mathcal F_1$: three consecutive weights]\label{ex:F1}
Let $3\le n\le m$.  Let $U'\le_{\fq}\fqm^2$ be a scattered subspace with $\dimq U'=n-1$, fix $\widetilde u\in U'\setminus\{0\}$ and $\lambda\in\fqm\setminus\fq$.  Since $U'\cap\langle\widetilde u\rangle_{\fqm}=\langle\widetilde u\rangle_{\fq}$, we have $\lambda\widetilde u\notin U'$, and we put
\[
  U=U'\oplus\langle\lambda\widetilde u\rangle_{\fq},
  \qquad \dimq U=n .
\]
Then:
\begin{itemize}
  \item $w_U(\langle\widetilde u\rangle_{\fqm})=2$, since $\widetilde u$ and $\lambda\widetilde u$ are $\fq$-independent;
  \item $w_U(P)\le w_{U'}(P)+1\le2$ for every point $P$, since $U'$ is an $\fq$-hyperplane of $U$;
  \item there are points of weight zero, by Lemma~\ref{lem:covering}, since $n\le m$;
  \item there are points of weight one.  Otherwise every point of $L_U$ would have weight two, and $|L_U|=(q^n-1)/(q^2-1)$; moreover every $2$-dimensional subspace $U\cap P$ would meet the hyperplane $U'$ nontrivially, so that $L_U=L_{U'}$ and $|L_U|=(q^{n-1}-1)/(q-1)$.  These two values coincide only for $n=2$.
\end{itemize}
Since $U'$ is scattered of rank $n-1\ge2$, it is not contained in a point, so $\langle U\rangle_{\fqm}=\fqm^2$.  Therefore the code $\cC$ associated with $U$ is a nondegenerate $[n,2]_{q^m/q}$ code with nonzero weights
\[
  n-2,\qquad n-1,\qquad n,
  \qquad g_1=g_2=1 .
\]
If $\nu=N_2(U)\ge1$, then \eqref{eq:vector-counting} gives $N_1(U)=\frac{q^n-1}{q-1}-(q+1)\nu$, hence
\[
  A_{n-2}(\cC)=(q^m-1)\nu,
  \qquad
  A_{n-1}(\cC)=(q^m-1)\Bigl(\frac{q^n-1}{q-1}-(q+1)\nu\Bigr).
\]
Thus, for every $q$, every $m\ge3$ and every $3\le n\le m$, there exists an antipodal three-weight $[n,2]_{q^m/q}$ code with no gaps.
\end{example}

The next family relies on the following explicit club.
\begin{example}[]\label{ex:m-2Club}
Let $m\ge3$, let $\lambda\in\fqm$ with $\fqm=\fq(\lambda)$, and let
\[
  U'=\bigl\{\bigl(a_0+a_1\lambda+\cdots+a_{m-2}\lambda^{m-2},\ a_{m-2}+a_{m-1}\lambda\bigr):a_0,\ldots,a_{m-1}\in\fq\bigr\}.
\]
Then $\dimq U'=m$, $\dimq\bigl(U'\cap\langle(1,0)\rangle_{\fqm}\bigr)=m-2$ and $w_{U'}(P)\le1$ for every point $P\ne\langle(1,0)\rangle_{\fqm}$.  In particular, for $m\ge4$, $L_{U'}$ is an $(m-2)$-club of rank $m$ with head $\langle(1,0)\rangle_{\fqm}$ (see \cite[Lemma 2.12]{DeBoeckVandeVoorde2016}).
\end{example}
\begin{example}[Family $\mathcal F_2$: weights $2$, $n-1$, $n$]\label{ex:F2}
Let $4\le n\le m$ and let $U'$ be as in Example~\ref{ex:m-2Club}.  Put $S_{m-2}=U'\cap\langle(1,0)\rangle_{\fqm}$, choose an $(n-2)$-dimensional $\fq$-subspace $S_{n-2}\le S_{m-2}$, and choose $u_1,u_2\in U'$ with $\langle u_1,u_2\rangle_{\fq}\cap S_{m-2}=\{0\}$ (for instance $u_1=(\lambda^{m-2},1)$ and $u_2=(0,\lambda)$).  Let
\[
  U_n=S_{n-2}\oplus\langle u_1,u_2\rangle_{\fq}.
\]
For instance, with the choices above and $S_{n-2}=\{(a_0+\cdots+a_{n-3}\lambda^{n-3},0) \,: \; a_i\in \fq\}$,
\[
  U_n=\bigl\{\bigl(a_0+a_1\lambda+\cdots+a_{n-3}\lambda^{n-3}+a_{m-2}\lambda^{m-2},\ a_{m-2}+a_{m-1}\lambda\bigr):a_i\in\fq\bigr\}.
\]
Then $\dimq U_n=n$ and $U_n\cap\langle(1,0)\rangle_{\fqm}=S_{n-2}$, because an element $s+\alpha u_1+\beta u_2$ of $U_n$ in $\langle(1,0)\rangle_{\fqm}$ lies in $S_{m-2}$, which forces $\alpha u_1+\beta u_2\in S_{m-2}$, hence $\alpha=\beta=0$.  Since $U_n\le U'$, every other point has weight at most one; the point $\langle u_1\rangle_{\fqm}$ has weight one and, as $n\le m$, there are points of weight zero.  So $L_{U_n}$ is an $(n-2)$-club of rank $n$, and the associated code $\cC$ is a nondegenerate $[n,2]_{q^m/q}$ code with nonzero weights
\[
  2,\qquad n-1,\qquad n,
  \qquad g_1=1,\quad g_2=n-3 .
\]
By \eqref{eq:vector-counting}, $N_{n-2}(U_n)=1$ and $N_1(U_n)=\frac{q^n-q^{n-2}}{q-1}=q^{n-2}(q+1)$, so
\[
  A_2(\cC)=q^m-1,
  \qquad
  A_{n-1}(\cC)=(q^m-1)q^{n-2}(q+1),
\]
and $A_n(\cC)=q^{2m}-1-(q^m-1)\bigl(1+q^{n-2}(q+1)\bigr)$.
For $n=3$ the construction still works, but the weights $2$ and $n-1$ coincide.  For $n\ge5$ the second gap is at least two, and it can be as large as $m-3$; however, unlike the first gap in Corollary~\ref{cor:antipodal-puncturing-three}, a large second gap does not impose any restriction on the parameters: the construction is valid for every $q$ and every $4\le n\le m$.
\end{example}

\begin{remark}\label{rem:clubs-general}
More generally, if $L_{U'}$ is an $i$-club of rank $N\le m$ in $\PG(1,q^m)$, with $2\le i\le N-1$, then by \eqref{eq:antipodal-point-weight} the code associated with $U'$ is a nondegenerate $[N,2]_{q^m/q}$ code with nonzero weights
\[
  N-i,\qquad N-1,\qquad N,
  \qquad g_1=1,\quad g_2=i-1,
\]
and, by \eqref{eq:vector-counting}, $A_{N-i}=q^m-1$ and $A_{N-1}=(q^m-1)\frac{q^N-q^i}{q-1}$.  The procedure of Example~\ref{ex:F2} produces new clubs from a given one: if $L_{U'}$ is an $i$-club of rank $m$ with head $P$, $S\le U'\cap P$ has dimension $j\ge2$ and $T\le U'$ has dimension $s\ge1$ with $T\cap P=\{0\}$ (so that $s\le m-i$), then $L_{S\oplus T}$ is a $j$-club of rank $j+s$.  For instance, starting from the $(m-1)$-club $\{(x,\Tr_{q^m/q}(x)):x\in\fqm\}$ one obtains $(n-1)$-clubs of rank $n$, hence codes with weights $1,n-1,n$, for every $3\le n\le m$.  For the existence of $i$-clubs of given rank in $\PG(1,q^m)$ see \cite{MannaertSantonastasoZullo2025} and the references therein.
\end{remark}

\begin{example}[Family $\mathcal F_3$: weights $a$, $2a-1$, $2a$]\label{ex:F3}
Let $m$ be even.  By \cite[Theorem~4.5]{NapolitanoEtAl2022}, for every $q$ there exists an $m$-dimensional $\fq$-subspace $U'\le\fqm^2$ of \emph{complementary weights}, that is,
\[
  w_{U'}(\langle(1,0)\rangle_{\fqm})=w_{U'}(\langle(0,1)\rangle_{\fqm})=\frac m2,
  \qquad
  w_{U'}(P)\le1\ \text{for every other point }P.
\]
Let $2\le a\le m/2$, and let $S_1\le U'\cap\langle(1,0)\rangle_{\fqm}$ and $S_2\le U'\cap\langle(0,1)\rangle_{\fqm}$ be $\fq$-subspaces with $\dimq S_1=\dimq S_2=a$.  Put
\[
  U=S_1\oplus S_2\ (\le U'),
  \qquad \dimq U=2a .
\]
Then $U\cap\langle(1,0)\rangle_{\fqm}=S_1$ and $U\cap\langle(0,1)\rangle_{\fqm}=S_2$ have dimension $a$, every point $\langle(s_1,s_2)\rangle_{\fqm}$ with $s_1s_2\ne0$, $(s_1,s_2)\in U$, has weight one, since $U\le U'$, and there are points of weight zero, since $2a\le m$.  By \eqref{eq:vector-counting}, $N_a(U)=2$ and $N_1(U)=(q^a-1)^2/(q-1)$.  Hence the code $\cC$ associated with $U$ is a nondegenerate $[2a,2]_{q^m/q}$ code with nonzero weights
\[
  a,\qquad 2a-1,\qquad 2a,
  \qquad g_1=1,\quad g_2=a-1,
\]
and
\[
  A_a(\cC)=2(q^m-1),
  \qquad
  A_{2a-1}(\cC)=(q^m-1)\frac{(q^a-1)^2}{q-1},
  \qquad
  A_{2a}(\cC)=q^{2m}-1-A_a(\cC)-A_{2a-1}(\cC).
\]
Observe that in this family the second largest weight $w=2a-1$ may be coprime with $m$ (for instance, $m=6$ and $a=3$ give the weights $3,5,6$), whereas, when $q\ge m$ and $g_1\ge2$, Corollary~\ref{cor:antipodal-puncturing-three} gives $t\mid\gcd(w,m)$ with $t\ge2$.
\end{example}

The last family realizes the situation described in Corollary~\ref{cor:antipodal-puncturing-three} with $t>n-w$, showing that the first gap can be large without $\cC$ being an expansion.

\begin{example}[Family $\mathcal F_4$: weights $\ell t-t-1$, $\ell t-t$, $\ell t-1$]\label{ex:F4}
Let $m=th$ with $1<t<m$, and let $2\le\ell\le h$.  Let $U_t$ be an $\mathbb F_{q^t}$-scattered $\mathbb F_{q^t}$-subspace of $\fqm^2=\mathbb F_{(q^t)^h}^2$ with $\dim_{\mathbb F_{q^t}}U_t=\ell$, so that $\dimq(U_t\cap P)\in\{0,t\}$ for every point $P$.  The $[\ell,2]_{q^m/q^t}$ code $\cC_t$ associated with $U_t$ has nonzero weights $\ell-1$ and $\ell$, so $\Exp_{q^t/q}(\cC_t)$ is an $[\ell t,2]_{q^m/q}$ code with nonzero weights $t(\ell-1)$ and $\ell t$.

Let $U\le_{\fq}U_t$ be an $\fq$-hyperplane of $U_t$, so that
\[
  n=\dimq U=\ell t-1\le m .
\]
Choosing an $\fq$-basis of $U_t$ whose first $n$ vectors form a basis of $U$, we see that the code $\cC$ associated with $U$ is obtained, up to equivalence, by deleting one coordinate (for instance the last one) from $\Exp_{q^t/q}(\cC_t)$.  Since $U$ is a hyperplane of $U_t$, every point of $L_{U_t}$ has weight $t-1$ or $t$ with respect to $U$, and the points outside $L_{U_t}$ have weight zero.  If $x$ and $y$ are the numbers of points of weight $t$ and $t-1$ with respect to $U$, then $x+y=|L_{U_t}|=\frac{q^{\ell t}-1}{q^t-1}$ and, by \eqref{eq:vector-counting}, $x(q^t-1)+y(q^{t-1}-1)=q^{\ell t-1}-1$; hence
\[
  x=\frac{q^{(\ell-1)t}-1}{q^t-1}\ge1,
  \qquad
  y=q^{(\ell-1)t}.
\]
Therefore $\cC$ is a nondegenerate $[\ell t-1,2]_{q^m/q}$ code with nonzero weights
\[
  d=\ell t-t-1=n-t,\qquad w=\ell t-t=n-t+1,\qquad n=\ell t-1,
  \qquad g_1=t-1,\quad g_2=1,
\]
and
\[
  A_{n-t}(\cC)=(q^m-1)\frac{q^{(\ell-1)t}-1}{q^t-1},
  \qquad
  A_{n-t+1}(\cC)=(q^m-1)\,q^{(\ell-1)t}.
\]
{If $t\geq 3$ (so that $g_1=t-1\geq 2$) and $q\geq m$, then} this family is an instance of Corollary~\ref{cor:antipodal-puncturing-three}.  Indeed, $\langle U\rangle_{\mathbb F_{q^t}}=U_t$, because $\langle U\rangle_{\mathbb F_{q^t}}\le U_t$ is an $\mathbb F_{q^t}$-subspace, so that its $\fq$-dimension is a multiple of $t$, and it contains $U$, of dimension $\ell t-1$.  The divisor $t$ satisfies $n-w=t-1<t\le\min\{w,m-w\}$ and $t\mid w$, exactly one coordinate is deleted, and $d=w-1>w-t$, so that $\cC_t$ has two weights and $A_w(\Exp_{q^t/q}(\cC_t))=A_d(\cC)+A_w(\cC)=(q^m-1)\frac{q^{w+t}-1}{q^t-1}$; the values of $A_d(\cC)$ and $A_w(\cC)$ above agree with the formulas of Corollary~\ref{cor:antipodal-puncturing-three}(iv).  Moreover, $\gcd(n-w,n-d)=\gcd(t-1,t)=1$, so $\cC$ is not the $\fq$-expansion of a code over any proper extension of $\fq$.  For $t=2$ we get three consecutive weights $2\ell-3,2\ell-2,2\ell-1$, as in Family~$\mathcal F_1$; {in this case $g_1=1$ and Corollary~\ref{cor:antipodal-puncturing-three} does not apply.  The smallest instance with $g_1\ge2$ is $m=6$, $t=3$, $\ell=2$: it gives a $[5,2]_{q^6/q}$ code with nonzero weights $2,3,5$ and first gap $2$, which is not an $\fq$-expansion from any proper extension of $\fq$ and which, for $q\ge7$, falls under Corollary~\ref{cor:antipodal-puncturing-three} with $t=3$.}
\end{example}

\begin{remark}\label{rem:gap-behaviour}
The following table summarizes the families above.
\begin{center}
\renewcommand{\arraystretch}{1.3}
\begin{tabular}{lllll}
\hline
Family & nonzero weights $d<w<n$ & $g_1=n-w$ & $g_2=w-d$ & parameters\\
\hline
$\mathcal F_1$ & $n-2,\ n-1,\ n$ & $1$ & $1$ & $3\le n\le m$\\
$\mathcal F_2$ & $2,\ n-1,\ n$ & $1$ & $n-3$ & $4\le n\le m$\\
$i$-clubs & $N-i,\ N-1,\ N$ & $1$ & $i-1$ & $2\le i\le N-1$, $N\le m$\\
$\mathcal F_3$ & $a,\ 2a-1,\ 2a$ & $1$ & $a-1$ & $m$ even, $2\le a\le m/2$\\
$\mathcal F_4$ & $\ell t-t-1,\ \ell t-t,\ \ell t-1$ & $t-1$ & $1$ & $t\mid m$, $2\le\ell\le m/t$\\
\hline
\end{tabular}
\end{center}
All these families exist for every $q$.  When $g_1=1$ (families $\mathcal F_1$, $\mathcal F_2$, $\mathcal F_3$ and clubs), Theorem~\ref{thm:antipodal-puncturing} does not apply, and indeed no arithmetic restriction appears: for every $q$ and every admissible length the second gap takes the values $1$ (Family $\mathcal F_1$), $n-3$ (Family $\mathcal F_2$) and $n-2$ (the $(n-1)$-clubs of Remark~\ref{rem:clubs-general}), and more generally the value $i-1$ whenever an $i$-club of rank $n$ exists; moreover, the weights need not share any common factor with $m$ (Family $\mathcal F_3$).  When $g_1\ge2$ and $q\ge m$, the situation is rigid: by Corollary~\ref{cor:antipodal-puncturing-three} there is a divisor $t$ of $m$ with $g_1\le t\le\min\{w,m-w\}$ dividing $w$, the residue $d\bmod t$ is either $0$ or at least $g_1$, and, if $t>g_1$, the weight distribution is determined by the parameters.  Family $\mathcal F_4$ shows that the case $t>g_1$ really occurs, with $t=g_1+1$, that then $\cC$ need not be an expansion, and that the first gap $g_1=t-1$ can be as large as $\frac m2-1$; in this family the second gap is equal to one.
\end{remark}
\endgroup

\section{{Simplex summands and field expansion}}\label{sec:codes}

In this section we state the consequences of Section~\ref{sec:linear-sets} for rank-metric codes.  Throughout, $\cC$ is a nondegenerate $[n,k]_{q^m/q}$ code with generator matrix $G$, so that every codeword is $xG$ for a unique $x\in\fqm^k$, and we put $N=km-n$.
By Proposition~\ref{prop:maximum-rank-weight}, $\cC$ contains codewords of rank $m$ if and only if $n\ge m$.  The key object is the subcode generated by the codewords of nonmaximum rank,
\[
  \cC_{<m}=\bigl\langle c\in\cC:0<\wt(c)<m\bigr\rangle_{\fqm},
  \qquad
  h(\cC)=\dimqm\cC_{<m}.
\]

\begin{lemma}\label{lem:code-dichotomy}
Let $\cC$ be a nondegenerate $[n,k]_{q^m/q}$ code with $h=h(\cC)\ge1$.  Then $\cC_{<m}$ contains a codeword of rank $m$ if and only if
\[
  n\ge(k-h+1)m.
\]
\end{lemma}

\begin{proof}
{If $h=0$, then $\cC_{<m}=\{0\}$ contains no codeword of rank $m$, and $n\le km<(k+1)m$; so we may assume $h\ge1$.}
Let $U$ be the system associated with $G$, and put $W=U^{\perp'}$ and $S=\langle W\rangle_{\fqm}$.  By \eqref{eq:dual-point-code-weight}, a nonzero vector $x$ satisfies $\wt(xG)<m$ if and only if $\langle x\rangle_{\fqm}\in L_W$.  The points of $L_W$ span $\PG(S,\fqm)$, and $x\mapsto xG$ is an $\fqm$-linear isomorphism from $\fqm^k$ onto $\cC$; hence
\begin{equation}
\label{eq: set-subcode-rk m}
  \cC_{<m}=\{xG:x\in S\}
  \qquad\text{and}\qquad
  \dimqm S=h.
\end{equation}
A codeword $xG$ with $x\in S\setminus\{0\}$ has rank $m$ if and only if $\langle x\rangle_{\fqm}\in\PG(S,\fqm)\setminus L_W$.  Therefore $\cC_{<m}$ contains a codeword of rank $m$ if and only if $L_W\ne\PG(S,\fqm)$, which by Lemma~\ref{lem:covering} is equivalent to $N\le(h-1)m$, that is, to $n\ge(k-h+1)m$.
\end{proof}

Accordingly, we distinguish two cases: either every nonzero codeword of $\cC_{<m}$ has rank smaller than $m$, equivalently $n<(k-h+1)m$ (this includes the case $h=0$), or $\cC_{<m}$ contains codewords of rank $m$, equivalently $n\ge(k-h+1)m$.  The first case is the subject of Section~\ref{subsec:direct-sum}, and the second that of Section~\ref{subsec:not-direct-sum}.  If $\cC$ is generated by its codewords of nonmaximum rank, then $h=k$ and the second case is simply $n\ge m$.

\subsection{Direct sum of a simplex code and a short code}\label{subsec:direct-sum}

We begin with the case in which $\cC_{<m}$ contains no codeword of rank $m$.
Geometrically, this means that the linear set $\Luperpp$ fills its span (see the proof of Lemma~\ref{lem:code-dichotomy}). We show that this happens precisely when $\cC$ splits off a simplex code, all of whose nonzero codewords have rank $m$, leaving a code of length smaller than $m$ that carries all the nonmaximum weights.
We first prove this in geometric terms (Theorem~\ref{thm:simplex-short-code}), and then express the condition through the parameters of $\cC$ (Lemma~\ref{lem: equiv-cond-direct-sum} and Corollary~\ref{cor:direct-sum-weights}).

\begin{theorem}\label{thm:simplex-short-code}
    Let $\cC$ be a nondegenerate $[n, k]_{q^m/q}$ rank-metric code and let $U$ be an associated $q$-system. Then there exists an $h \in \{0,\ldots, k\}$ such that
    \[\cC \cong \tilde{\cC} \times \cC'=\{(\tilde c,c'):\tilde c\in\tilde\cC,\ c'\in\cC'\}\le\fqm^{(k-h)m} \times \fqm^{\varepsilon}\] where $\tilde{\cC}$ is an $[(k-h)m, k-h]_{q^m/q}$ simplex code and $\cC'$ is an $[\varepsilon,h]_{q^m/q} $ rank-metric code of short length $0 \le \varepsilon \le m-1$ if and only if
    $$L_{U^{\perp'}} = \PG(S, \fqm)$$ for some $\fqm$-subspace $S \subseteq \fqmk$ with $\dim_{\fqm}(S) = h$. For $h=0$ \textup(resp.\ $h=k$\textup), $\cC'$ \textup(resp.\ $\tilde\cC$\textup) is the zero code of length $0$.  Moreover, in this case  $n=(k-h)m+\varepsilon$ and the code $\cC'$ is nondegenerate; furthermore, $S=\langle U^{\perp'}\rangle_{\fqm}$ and $h=h(\cC)$.
\end{theorem}
\begin{proof} We first observe that the condition on $\Luperpp$ does not depend on the choice of the $q$-system $U$ associated with $\cC$.  Any other such system is $\varphi(U)$, with $\varphi\in\GL(k,q^m)$; if $\psi\in\GL(k,q^m)$ is the adjoint of $\varphi^{-1}$ with respect to $\sigma$, then $\varphi(U)^{\perp'}=\psi(\Uperpp)$, and $\psi$, being $\fqm$-linear, maps $\PG(S,\fqm)$ onto $\PG(\psi(S),\fqm)$.\\
  ($\Leftarrow$)   Suppose that $S$ is an $h$-dimensional $\fqm$-subspace of $\fqmk$ with $h\in \{0,\ldots, k\}$ such that $$\Luperpp = \PG(S,\fqm).$$
  If $h=0$, we get that $U = \fqmk$ and therefore $\cC$ is a simplex code. Let $h \ge 1$. Then we have that $\Uperpp$ is an $\fq$-subspace of $S$ and since $\dimq(U) = n$, we have that \begin{equation}
    \label{eq: dim-Uperpp1}
        km - n = \dimq(\Uperpp) \le \dimq(S)= hm
    \end{equation} and by Lemma \ref{lem:covering} 
\begin{equation}
    \label{eq: dim-Uperpp2}
        km - n = \dimq(\Uperpp) \ge (h-1)m+1. 
    \end{equation} 
Combining \eqref{eq: dim-Uperpp1} and \eqref{eq: dim-Uperpp2}, we get the following: 
\[ (k-h)m  \le n \le (k-h)m + m - 1.\] 
So $n = (k-h)m + \varepsilon$ for $0 \le \varepsilon \le m-1$.
The subspace $S$ is $\fqm$-linear and $\Uperpp \leqslant_{\fq} S$. Hence $ S^\perp = S^{\perp'} \leqslant_{\fq}(\Uperpp)^{\perp'} = U$ , and $\dimqm(\Sp)=k-h$.  We want to construct an $\varepsilon$-dimensional subspace $U'$ such that \begin{equation}
\label{eq: sum-decomposition-subspaces}
    U = S^\perp \oplus U' \text{ and } S^\perp \cap \langle U' \rangle_{\fqm} = \{0\}.
\end{equation}
\\
We start the construction by taking any $\fq$-subspace $\Us$ of $U$ such that 
\[U=S^\perp \oplus U^*.\]
Then we get
$$\langle U \rangle_{\fqm} = S^\perp + \Usqm.$$
Let $g = \dimqm(\Usqm)$. We have $\dimqm(\Sp) = k-h$ and since $\cC$ is nondegenerate, $\langle U \rangle_{\fqm} = \fqmk$. Therefore $g \ge h$ and $\dimqm(S^\perp \cap \langle \Us \rangle_{\fqm} ) = g-h.$ We define $H \leqslant_{\fqm} \Usqm$ of $\dimqm(H) = h$ such that $$H \cap (S^\perp \cap \langle \Us\rangle_{\fqm}) = \{0\}.$$ 
%\textcolor{blue}{(such an $H$ exists, since $S^\perp\cap\Usqm$ has codimension $h$ in $\Usqm$; as $H\subseteq\Usqm$, this means $H\cap\Sp=\{0\}$)}
Over $\fq$, we define the subspace $$\overline{U} = U \cap \Usqm,$$ since $U+\Usqm\le\fqmk$,  we have that $\dimq(U) +  \dimq(\Usqm)  - \dimq(\overline{U}) \le km$. In particular, we compute $$\dimq(\overline{U}) \ge (k-h)m+\varepsilon + gm - km = (g-h)m + \varepsilon.$$ 
We define the desired subspace $U'$ as the intersection $$U' = H \cap \overline{U}.$$ 
Since $U' \subseteq H$, also $\Upqm \subseteq H$ and therefore $\Upqm \cap \Sp = \{0\}$. Since $U' \subseteq \overline{U} \subseteq U$ and $\Sp \subseteq U$,
\begin{align*}
    & \dimq(U') + \dimq(\Sp) = \dimq(U' + \Sp)\le \dimq(U) \\ \Leftrightarrow & \dimq(U') \le (k-h)m + \varepsilon - (k-h)m = \varepsilon.
\end{align*} 
Furthermore, since $H+\overline U\subseteq\Usqm$, we compute \begin{align*}
    &gm \ge \dimq(H+ \overline U) \ge hm + (g-h)m+\varepsilon - \dimq(H\cap \overline{U})\\ \Leftrightarrow \,\, &gm \ge \dimq(H+ \overline U) \ge gm + \varepsilon - \dimq(H\cap \overline{U}).
\end{align*}
Hence $\dimq(U') = \dimq(H \cap \overline{U}) = \varepsilon$ and we conclude that $U = \Sp \oplus U'.$ Moreover, $\fqmk=\langle U\rangle_{\fqm}=\Sp+\Upqm$, and $\Upqm\subseteq H$ with $\dimqm(H)=h$; comparing dimensions, $\Upqm=H$. Hence
 $$\fqmk = \Sp \oplus \Upqm,$$ and we may choose an $\fqm$-basis $B = (v_1, \hdots , v_{k-h}, v_{k-h+1},\hdots, v_k)$ of $\fqmk$ such that $$\Sp = \langle v_i : 1 \le i \le k-h\rangle_{\fqm}  \text{ and } \Upqm = \langle v_i : k-h+1 \le i \le k\rangle_{\fqm}.$$
Let $\phi \in \GL(k,q^m)$ be the linear transformation mapping $B$ to the canonical basis $(e_1, \hdots, e_k)$ of $\fqmk$. Since $\phi$ induces an equivalence of codes, we have that $$\phi(U) = \phi(\Sp) \oplus \phi(U')$$ is the system associated to a code equivalent to $\cC$, with $\phi(\Sp)=\fqm^{k-h}\times\{0\}$ and $\phi(U')\subseteq\{0\}\times\fqm^h$. Now we choose corresponding $\fq$-bases 
$B_{\phi(\Sp)} = \{ (\overline{z}_1,\overline{0}), \hdots , (\overline{z}_{(k-h)m},\overline{0}) \}$ of   $\phi(\Sp)$ and $B_{\phi(U')} = \{ (\overline{0},\overline{y}_1), \hdots , (\overline{0},\overline{y}_\varepsilon) \} \text{ of } \phi(U')$, where $\overline z_i\in\fqm^{k-h}$ and $\overline y_j\in\fqm^h$.
With respect to these bases, the generator matrix of the code associated to $U$ has the block form

$$ G = \begin{pmatrix}
    z_1^{(1)} & \hdots & z_{(k-h)m}^{(1)} & 0 & \hdots & 0\\ \vdots & & \vdots & \vdots & & \vdots \\ z_1^{(k-h)} & \hdots & z_{(k-h)m}^{(k-h)} & 0 & \hdots & 0 \\ 0 & \hdots & 0 & y_1^{(1)} & \hdots & y_\varepsilon^{(1)} \\  \vdots & & \vdots & \vdots & & \vdots\\ 0 & \hdots & 0 & y_1^{(h)} & \hdots & y_\varepsilon^{(h)}
\end{pmatrix} \in \fqm^{k\times n}.$$
The system $\phi(\Sp)=\fqm^{k-h}\times\{0\}$ is $\fqm$-linear of dimension $k-h$. Hence the code $\tilde\cC$ generated by the block $Z=(\overline z_1\ \cdots\ \overline z_{(k-h)m})$ is a nondegenerate $[(k-h)m, k-h]_{q^m/q}$ simplex code. The code $\cC'$ generated by the block $Y=(\overline y_1\ \cdots\ \overline y_\varepsilon)$ is a nondegenerate $[\varepsilon,h]_{q^m/q}$ code, since $\overline y_1,\ldots,\overline y_\varepsilon$ are $\fq$-linearly independent and span $\fqm^h$ over $\fqm$ (because $\phi(\Upqm)=\{0\}\times\fqm^h$).  As $G$ is block diagonal, the code generated by $G$ is $\tilde\cC\times\cC'$, and therefore $\cC\cong\tilde\cC\times\cC'$.  For $h=k$ the same argument works with $\Sp=\{0\}$ and $Z$ empty, and gives $\cC\cong\cC'$. This proves the first implication. \\
\\
{($\Rightarrow$)} 
Conversely, suppose that there exists an $h \in \{0,\hdots, {k}\}$ such that
    \[\cC \cong \tilde{\cC} \times \cC' \subseteq \fqm^{(k-h)m} \times \fqm^\varepsilon, \] with $0 \le \varepsilon \le m-1$ and $\dimqm (\tilde{\cC }) = k-h $ and $\dimqm(\cC') = h$.
   By the first paragraph, we may compute $U$ from any generator matrix of $\tilde\cC\times\cC'$.  Write an $\fqm$-basis of $\tilde\cC\times\cC'$ where $(c_1, 0) , \hdots , (c_{k-h},0) \in \fqm^{(k-h)m} \times \{0\}$ and $(0, c_{k-h+1}) , \hdots , (0,c_k) \in \{0\} \times \fqm^\varepsilon$, and let $\tilde G$ and $G'$ be the matrices with rows $c_1,\ldots,c_{k-h}$ and $c_{k-h+1},\ldots,c_k$, respectively.  Then $\mathrm{diag}(\tilde G,G')$ is a generator matrix of $\tilde\cC\times\cC'$, and the associated system is the $\fq$-span of its columns:
     $$U = \bigl(\tilde{U}\times\{0\}\bigr) \oplus \bigl(\{0\}\times U'\bigr),$$
    where $\tilde U\le\fqm^{k-h}$ and $U'\le\fqm^{h}$ are the $\fq$-spans of the columns of $\tilde G$ and $G'$.  Since $\cC$ is nondegenerate, these $n$ columns are $\fq$-linearly independent, so $n=(k-h)m+\varepsilon$ and $\cC'$ is nondegenerate.
     Since $\tilde{\cC}$ is a simplex code, the associated system $\tilde{U} = \fqm^{k-h}$ is $\fqm$-linear. Let $$S = \{0\}^{k-h} \times \fqm^h.$$ Then $S$ is an $h$-dimensional $\fqm$-subspace such that $U^{\perp'} \subseteq {(\tilde{U}\times\{0\})^{\perp'} = (\tilde{U}\times\{0\})^\perp} = S$, and
     $$\dimq(U^{\perp'}) = km - n = km - (k-h)m - \varepsilon = hm - \varepsilon \ge hm - (m-1) = (h-1)m+1.$$
     So ${\dimq(\Uperpp)} > (h-1)m$ and since $\Uperpp\le_{\fq}S$, if $h\ge1$ Lemma~\ref{lem:covering} gives $L_{\Uperpp} = \PG(S,\fqm)$.  If $h=0$, then $\Uperpp\subseteq S=\{0\}$ and $\Luperpp=\varnothing=\PG(S,\fqm)$.
    
     \end{proof}

\begin{lemma}
    \label{lem: equiv-cond-direct-sum}
    Let $\cC$ be a nondegenerate $[n,k]_{q^m/q}$ code, let $U$ be an associated system and put $h=h(\cC)$.  Then
\begin{equation}
\label{eq: h lower bound}
  h\ge k-\left\lfloor\frac nm\right\rfloor,
\end{equation}
and the following are equivalent:
\begin{enumerate}[label=\textup{(\alph*)}]
  \item $\cC_{<m}$ contains no codeword of rank $m$;
  \item $n\le(k-h+1)m-1$;
  \item $h=k-\lfloor n/m\rfloor$;
  \item $L_{U^{\perp'}}=\PG(\langle U^{\perp'}\rangle_{\fqm},\fqm)$.
\end{enumerate}
If $h=k$, these conditions mean that $n\le m-1$, and then $\cC$ has no codeword of rank $m$.
\end{lemma}

\begin{proof}
Let $S=\langle \Uperpp \rangle_{\fqm}$. We use the setting of Lemma \ref{lem:code-dichotomy}. We write $\cC_{<m} = \{xG : x \in S\}$ and $\dimqm S = h$. Since $\Uperpp$ is a subspace of $S$, we have $km-n \le hm$. Hence, \[
  h\ge k-\left\lfloor\frac nm\right\rfloor.
\]

 (a) $\Leftrightarrow $ (b): see Lemma \ref{lem:code-dichotomy}.
 
    (b) $\Leftrightarrow $ (c):  $n \le (k-h+1)m -1 \Leftrightarrow n < (k -h + 1)m \Leftrightarrow \frac nm < k-h+1 \Leftrightarrow \left\lfloor \frac{n}{m} \right\rfloor \le k-h$. 
    
    Together with \eqref{eq: h lower bound}, $n \le (k-h+1)m -1 \Leftrightarrow \left\lfloor \frac{n}{m} \right\rfloor = k-h $.

    (b) $\Leftrightarrow$ (d): $n < (k-h+1)m \Leftrightarrow (h-1)m < km - n = \dimq(\Uperpp) \Leftrightarrow L_{\Uperpp} = \PG(S,\fqm)$, 
    
    see Lemma \ref{lem:covering}. 
    
    \noindent If $h=k$; (b) $ \Rightarrow n\le m-1$. By Proposition \ref{prop:maximum-rank-weight}, we conclude that $\cC$ has no codeword of rank $m$. 
   
\end{proof}

\begin{corollary}\label{cor:direct-sum-weights}

Let $\cC$ be a nondegenerate $[n,k]_{q^m/q}$ code, let $U$ be an associated system and put $h=h(\cC)$. Suppose that the equivalent conditions in Lemma \ref{lem: equiv-cond-direct-sum} hold, then, up to equivalence, $$\cC=\tilde\cC\times\cC',$$ where $\tilde\cC$ is a $[(k-h)m,k-h]_{q^m/q}$ simplex code (the zero code of length $0$ if $h=k$) and $\cC'$ is a nondegenerate $[\varepsilon,h]_{q^m/q}$ code with $\varepsilon=n-(k-h)m\le m-1$. Furthermore, $\cC_{<m}=\{0\}\times\cC'$, and the nonzero weights of $\cC$ are the nonzero weights of $\cC'$, together with $m$ if $h\le k-1$; the weights of $\cC'$ are at most $\varepsilon\le m-1$, and
\[
  A_m(\cC)=q^{km}-q^{hm},
  \qquad
  A_i(\cC)=A_i(\cC')\quad\text{for } 1\le i\le m-1 .
\]

% and for every $(\tilde c,c')\in\tilde\cC\times\cC'$
% \[
%   \wt(\tilde c,c')=
%   \begin{cases}
%     m & \text{if } \tilde c\ne0,\\
%     \wt(c') & \text{if } \tilde c=0.
%   \end{cases}
% \]
 
\end{corollary}

\begin{proof}

    If the conditions of Lemma \ref{lem: equiv-cond-direct-sum} hold, then, by Theorem \ref{thm:simplex-short-code}, up to equivalence, $\cC=\tilde\cC\times\cC'$, where $\tilde\cC$ is a $[(k-h)m,k-h]_{q^m/q}$ simplex code \textup(the zero code of length $0$ if $h=k$\textup) and $\cC'$ is a nondegenerate $[\varepsilon,h]_{q^m/q}$ code with $\varepsilon=n-(k-h)m\le m-1$.  
     Since every nonzero codeword of a simplex code has rank weight $m$, we get that for every $(\tilde c,c')\in\tilde\cC\times\cC'$
\[
  \wt(\tilde c,c')=
  \begin{cases}
    m & \text{if } \tilde c\ne0,\\
    \wt(c') & \text{if } \tilde c=0.
  \end{cases}
\]
Every codeword of $\cC'$ has rank weight at most $\varepsilon \le m-1$. Therefore, $\cC_{<m} = \{0\} \times \cC'$ and $A_i(\cC)=A_i(\cC')\quad\text{for } 1\le i\le m-1 $. The nonzero weights of $\cC$ are the nonzero weights of $\cC'$, together with $m$ if $h\le k-1$, since then $\dimqm(\tilde{\cC}) = (k-h) > 0$. Counting the codewords gives the relation $ \sum\limits_{i=1}^m A_i(\cC)= A_m(\cC) + \sum\limits_{i=1}^{m-1} A_i(\cC').$ Hence $$A_m(\cC) = \sum\limits_{i=1}^m A_i(\cC)- \sum\limits_{i=1}^{m-1} A_i(\cC') = q^{km} - 1 - (q^{hm}-1) = q^{km} - q^{hm}. $$ 
\end{proof}

\subsection{Codes that are not a direct sum of a simplex code and a short code}\label{subsec:not-direct-sum}

Throughout this subsection, the nonzero rank weights of $\cC$ are $d_1<\cdots<d_r<m$ together with $m$, we put
\[
  h=h(\cC),
  \qquad
  a=m-d_r,
  \qquad
  N=km-n,
\]
and we assume:
\begin{enumerate}[label=\textup{(H\arabic*)},ref=\textup{(H\arabic*)}]
  \item\label{H1} $q\ge m$;
  \item\label{H2} $h\ge2$;
  \item\label{H3} $\cC_{<m}$ contains a codeword of rank $m$; equivalently, $n\ge(k-h+1)m$, by Lemma~\ref{lem:code-dichotomy};
  \item\label{H4} $a\ge2$, that is, $d_r\le m-2$.
\end{enumerate}
Note that \ref{H2} forces $r\ge1$, and \ref{H3} forces the weight $m$ to occur.

\begin{remark}\label{rem:code-hypotheses}
Assumption \ref{H2} holds as soon as $r\ge2$.  Indeed, if $h\le1$, then $\cC_{<m}$ is generated by a single codeword, and all nonzero $\fqm$-multiples of a codeword have the same rank, so at most one weight smaller than $m$ occurs.  If $h=k$, assumption \ref{H3} is equivalent to $n\ge m$; and if $n>(k-1)m$, then \ref{H2} implies \ref{H3}, since $(k-h+1)m\le(k-1)m$.  Finally, if \ref{H4} fails, then $d_r=m-1$, the greatest common divisor of the deficiencies $m-d_i$ equals one, and no reduction to a larger field can be expected.
\end{remark}

The first result shows that deleting a controlled number of coordinates turns $\cC$ into an expansion from a larger field, with an explicit effect on the ranks of the codewords.

\begin{theorem}\label{thm:larger-field-puncturing}
Assume \ref{H1}--\ref{H4}.  Then there exists a divisor $t$ of $m$ with
\[
  a\le t\le N-a
  \qquad\text{and}\qquad
  t\mid n-d_r,
\]
such that, after replacing $\cC$ by an equivalent code, the map
\[
  \pi:\cC\longrightarrow\fqm^{\,n-t+a}
\]
deleting the last $t-a$ coordinates has the following properties.
\begin{enumerate}[label=\textup{(\roman*)}]
\item The punctured code $\pi(\cC)$ is, up to equivalence, the $\fq$-expansion of a $\bigl[\tfrac{n-t+a}{t},k'\bigr]_{q^m/q^t}$ code, where $k'=\dimqm\pi(\cC)$.
\item For every $c\in\cC$,
\[
  \wt(\pi(c))=t\left\lfloor\frac{\wt(c)}{t}\right\rfloor.
\]
Equivalently, since $t\mid m$, the deficiency $m-\wt(\pi(c))$ is the least multiple of $t$ which is at least $m-\wt(c)$.  In particular, $\wt(\pi(c))=m$ if and only if $\wt(c)=m$.
\end{enumerate}
\end{theorem}

\begin{proof}
Let $U$ be the system associated with $G$, and put $W=U^{\perp'}$ and $S=\langle W\rangle_{\fqm}$.  As in the proof of Lemma~\ref{lem:code-dichotomy}, $\dimqm S=h\ge2$, and \ref{H3} gives $N=\dimq W\le(h-1)m$.  By \eqref{eq:dual-point-code-weight}, the positive point weights of $L_W$ are $m-d_r<\cdots<m-d_1$, whose minimum is $a\ge2$.  By \ref{H1}, Theorem~\ref{thm:larger-field-representation} applies to $W$: there are a divisor $t\mid m$ and $W_t=\langle W\rangle_{\mathbb F_{q^t}}$ such that
\[
  L_{W_t}=L_W,
  \qquad
  \dimq W_t=N+t-a,
  \qquad
  a\le t\le N-a,
  \qquad
  t\mid N-a.
\]
Since $N-a=(k-1)m-n+d_r$ and $t\mid m$, also $t\mid n-d_r$.

Put $U_t=W_t^{\perp'}$.  Since $W\subseteq W_t$, we have $U_t\subseteq U$ and
\[
  \dimq U_t=km-\dimq W_t=n-t+a;
\]
by Lemma~\ref{lem:perp-linear}, $U_t$ is $\mathbb F_{q^t}$-linear.  Choose an $\fq$-basis of $U$ whose first $n-t+a$ vectors form a basis of $U_t$.  The matrix having these vectors as columns is $GA$ for some $A\in\GL(n,q)$, and it generates a code equivalent to $\cC$; we replace $G$ by $GA$.  Then $\pi(xG)=xG_t$, where $G_t$ consists of the first $n-t+a$ columns of $G$, which form an $\fq$-basis of $U_t$.

(i) The code $\pi(\cC)$ is generated by the rows of $G_t$, whose rank is $k'=\dimqm\langle U_t\rangle_{\fqm}$.  If $k'<k$, we identify $\langle U_t\rangle_{\fqm}$ with $\fqm^{k'}$, so that $U_t$ becomes an $[n-t+a,k']_{q^m/q}$ system to which Definition~\ref{def:expansion} applies. Choose an $\mathbb F_{q^t}$-basis of $U_t$ and expand it with respect to a fixed $\fq$-basis of $\mathbb F_{q^t}$, as in Definition~\ref{def:expansion}.  The resulting vectors form another $\fq$-basis of $U_t$, so the matrix having them as columns generates a code equivalent to $\pi(\cC)$.  Hence $\pi(\cC)$ is, up to equivalence, an $\fq$-expansion from $\mathbb F_{q^t}$.

(ii) The case $c=0$ is trivial.  Let $x\in\fqm^k\setminus\{0\}$, $P=\langle x\rangle_{\fqm}$ and $c=xG$.  Applying \eqref{eq:dual-point-code-weight} to $U$ and to $U_t$ gives
\[
  \wt(c)=m-w_W(P),
  \qquad
  \wt(\pi(c))=m-w_{W_t}(P).
\]
If $\wt(c)=m$, then $P\notin L_W=L_{W_t}$, and hence $\wt(\pi(c))=m$.  Otherwise $w_W(P)=m-\wt(c)>0$.  Since $(W_t\cap P)/(W\cap P)$ embeds into $W_t/W$, which has dimension $t-a$,
\[
  m-\wt(c)\le w_{W_t}(P)\le m-\wt(c)+t-a.
\]
As $W_t$ is $\mathbb F_{q^t}$-linear, $w_{W_t}(P)$ is a multiple of $t$; since the interval above has length $t-a<t$, it contains at most one multiple of $t$.  The least multiple of $t$ which is at least $m-\wt(c)$ lies between $m-\wt(c)$ and $w_{W_t}(P)$, hence in the interval, and therefore equals $w_{W_t}(P)$.  Thus
\[
  \wt(\pi(c))
  =m-t\left\lceil\frac{m-\wt(c)}{t}\right\rceil
  =t\left\lfloor\frac{\wt(c)}{t}\right\rfloor,
\]
where the last equality holds because $t\mid m$.
\end{proof}

\begin{remark}\label{rem:two-directions}
Theorems~\ref{thm:antipodal-puncturing} and~\ref{thm:larger-field-puncturing} go in opposite directions.  In Theorem~\ref{thm:antipodal-puncturing}, Theorem~\ref{thm:larger-field-representation} is applied to the system $U$ itself; the subspace $\langle U\rangle_{\mathbb F_{q^t}}$ contains $U$, and $\cC$ is obtained by deleting coordinates \emph{from} an expansion.  In Theorem~\ref{thm:larger-field-puncturing}, it is applied to $U^{\perp'}$; enlarging $U^{\perp'}$ to $W_t$ reduces $U$ to $U_t=W_t^{\perp'}$, so that deleting coordinates \emph{of} $\cC$ produces an expansion.  For two-dimensional codes of length $n<m$ only the first approach is available, since $\dimq U^{\perp'}=2m-n>m$ and hence $L_{U^{\perp'}}=\PG(1,q^m)$ by Lemma~\ref{lem:covering}.  For $n=m$ and $k=2$ both apply, and both give $t=a$: for Theorem~\ref{thm:antipodal-puncturing} see Remark~\ref{rem:antipodal-consequences}(c), and in Theorem~\ref{thm:larger-field-puncturing} one has $t\mid n-d_r=a$.
\end{remark}

The next proposition describes the punctured code explicitly.

\begin{proposition}\label{prop:weight-transformation}
Assume \ref{H1}--\ref{H4}, and let $t$ and $\pi$ be as in Theorem~\ref{thm:larger-field-puncturing}.
\begin{enumerate}[label=\textup{(\roman*)}]
\item The map $\pi$ is injective, so that $\pi(\cC)$ is a nondegenerate $[n-t+a,k]_{q^m/q}$ code, if and only if $d_1\ge t$.

\item For every $i=1,\ldots,r$,
\[
  d_i\bmod t\le t-a.
\]
Moreover, if $a<t$, then $d_r\bmod t=t-a$.

\item Suppose that $d_1\ge t$.  Then the nonzero weights of $\pi(\cC)$ are $m$ together with the distinct values among
\[
  d_i^{(t)}=t\left\lfloor\frac{d_i}{t}\right\rfloor,
  \qquad i=1,\ldots,r;
\]
they are all multiples of $t$, and $d_r^{(t)}=m-t$.  Moreover,
\[
  A_m(\pi(\cC))=A_m(\cC),
\]
where $A_j$ denotes the number of codewords of rank $j$, and for every $j<m$,
\[
  A_j(\pi(\cC))
  =
  \sum_{\substack{1\le i\le r\\ d_i^{(t)}=j}}
      A_{d_i}(\cC).
\]
\end{enumerate}
\end{proposition}

\begin{proof}
By Theorem~\ref{thm:larger-field-puncturing}(ii), $\pi(c)=0$ if and only if $\wt(c)<t$.  Hence $\pi$ is injective if and only if $d_1\ge t$.  The columns of the punctured generator matrix are $\fq$-linearly independent, being columns of $G$; so, when $\pi$ is injective, $\pi(\cC)$ is a nondegenerate code of dimension $k$.  This proves~(i).

Deleting $t-a$ coordinates of a vector decreases its rank by at most $t-a$.  Hence, for a codeword $c$ of rank $d_i$, Theorem~\ref{thm:larger-field-puncturing}(ii) gives
\[
  d_i\bmod t
  =d_i-t\left\lfloor\frac{d_i}{t}\right\rfloor
  =\wt(c)-\wt(\pi(c))
  \le t-a.
\]
If $a<t$, then $d_r=m-a\equiv t-a\pmod t$, because $t\mid m$ and $0<t-a<t$.  This proves~(ii).

Under the assumption of~(iii), $\pi$ is a bijection from $\cC$ onto $\pi(\cC)$ which, by Theorem~\ref{thm:larger-field-puncturing}(ii), maps codewords of rank $d_i$ to codewords of rank $d_i^{(t)}$ and codewords of rank $m$ to codewords of rank $m$.  All these values are multiples of $t$, and $d_r^{(t)}=t\lfloor(m-a)/t\rfloor=m-t$ because $1\le a\le t$ and $t\mid m$.  The formulas for the multiplicities follow.
\end{proof}

\begin{remark}\label{rem:residue-condition}
Puncturing may decrease the number of distinct nonzero weights, since different weights may have the same largest multiple of $t$ below them.  Proposition~\ref{prop:weight-transformation}(ii) is a genuine restriction on the weight spectrum of codes satisfying \ref{H1}--\ref{H4}.  For instance, if $\cC$ has exactly three nonzero weights $d<w<m$, then there is a divisor $t$ of $m$ with
\[
  m-w\le t,
  \qquad
  t\mid n-w,
  \qquad
  d\bmod t\le t-(m-w).
\]
\end{remark}

\begin{example}\label{ex:weight-collapse}
Suppose that $m=12$, that a code $\cC$ satisfying \ref{H1}--\ref{H4} has nonzero rank weights
\[
  5,\ 6,\ 8,\ 10,\ 12,
\]
and that $t=4$ in Theorem~\ref{thm:larger-field-puncturing}.  Then $a=12-10=2$, so $t-a=2$ coordinates are deleted.  The residues of the nonmaximum weights modulo $4$ are $1,2,0,2$, all at most $t-a=2$, as required by Proposition~\ref{prop:weight-transformation}(ii).  The ranks transform as
\[
  5\longmapsto4,\qquad
  6\longmapsto4,\qquad
  8\longmapsto8,\qquad
  10\longmapsto8,\qquad
  12\longmapsto12.
\]
Since $d_1=5\ge t$, the punctured code is a nondegenerate $[n-2,k]_{q^{12}/q}$ code with nonzero weights $4,8,12$, and
\[
  A_4(\pi(\cC))=A_5(\cC)+A_6(\cC),
  \qquad
  A_8(\pi(\cC))=A_8(\cC)+A_{10}(\cC),
  \qquad
  A_{12}(\pi(\cC))=A_{12}(\cC).
\]
Thus four distinct nonmaximum weights collapse to two.  By contrast, $7$ cannot be a weight of such a code: $7\bmod4=3>t-a$.  Equivalently, a codeword of rank $7$ would have to lose three units of rank when only two coordinates are deleted.
\end{example}

We now bring in the divisor-gap condition of Definition~\ref{def:gap}.

\begin{corollary}\label{cor:code-rigidity}
Assume \ref{H1}--\ref{H4}, and put $e=\gcd(m-d_1,\ldots,m-d_r)$.
If $\cD_m(e,N)=\varnothing$, then $e=m-d_r$,
and $\cC$ is, up to equivalence, the $\fq$-expansion of an $[n/e,k]_{q^m/q^e}$ rank-metric code.  In particular,
\[
    e\mid m,
    \qquad
    e\mid n,
    \qquad
    e\mid d_i\ \text{ for every } i.
\]
Moreover, every divisor $t$ as in Theorem~\ref{thm:larger-field-puncturing} equals $e$, so that no coordinate is deleted.
\end{corollary}

\begin{proof}
Let $t$ be as in Theorem~\ref{thm:larger-field-puncturing}.  Since $e$ divides $m-d_r=a$,
\[
  e\le a\le t\le N-a\le N-e.
\]
As $t\mid m$ and $\cD_m(e,N)=\varnothing$, we get $t=e$, and then $e\le a\le t$ forces $a=t=e$.  Thus $t-a=0$, the map $\pi$ is the identity, and Theorem~\ref{thm:larger-field-puncturing}(i) shows that $\cC$ is, up to equivalence, the $\fq$-expansion of an $[n/e,k]_{q^m/q^e}$ code; in particular $e\mid n$, and $e\mid m$ because $t\mid m$.  Since expansion multiplies all rank weights by $e$, every $d_i$ is divisible by $e$.
\end{proof}

\begin{corollary}\label{cor:three-weight-rigidity}
Let $\cC$ be a nondegenerate $[n,k]_{q^m/q}$ rank-metric code with nonzero weights $d<w<m$ satisfying \ref{H1}, \ref{H3} and \ref{H4}; assumption \ref{H2} is then automatic by Remark~\ref{rem:code-hypotheses}.  Put $e=\gcd(m-w,m-d)$.  If $\cD_m(e,km-n)=\varnothing$, then $e=m-w$.  Consequently,
\[
    m-w\mid m-d,\qquad
    m-w\mid m,\qquad
    m-w\mid n,
\]
and $\cC$ is, up to equivalence, the $\fq$-expansion of an $\left[\frac{n}{m-w},k\right]_{q^m/q^{m-w}}$ rank-metric code whose nonzero weights are
\[
  \frac d{m-w},\qquad \frac w{m-w},\qquad \frac m{m-w}.
\]
\end{corollary}

\begin{proof}
This is Corollary~\ref{cor:code-rigidity} with $r=2$, $d_1=d$ and $d_2=w$.  Since expansion multiplies all rank weights by $e=m-w$, the weights of the shorter code are obtained by dividing those of $\cC$ by $m-w$.
\end{proof}

Corollary~\ref{cor:code-rigidity} requires the divisor-gap condition.  There is, however, an important case in which no such assumption is needed, namely when the length of the code is a multiple of $m$.

\begin{theorem}\label{thm:length-multiple-m}
Let $\cC$ be a nondegenerate $[n,k]_{q^m/q}$ rank-metric code whose nonzero rank weights are $d_1<\cdots<d_r<m$, with $r\ge1$, together with $m$.  Assume that $m\mid n$, and put $a=m-d_r$.  Then $a\mid m$, $a\mid n$, and $\cC$ is, up to equivalence, the $\fq$-expansion of an $[n/a,k]_{q^m/q^a}$ rank-metric code.  In particular, every rank weight of $\cC$ is divisible by $a$.
\end{theorem}

\begin{proof}
Let $U$ be an associated system and put $W=U^{\perp'}$.  Since $\dimq W=km-n$ and $m\mid n$, the rank of $W$ is a multiple of $m$.  By \eqref{eq:dual-point-code-weight}, the positive point weights of $L_W$ are $m-d_r<\cdots<m-d_1$, whose minimum is $a$.  Proposition~\ref{prop:multiple-m} yields $a\mid m$ and shows that $W$ is $\mathbb F_{q^a}$-linear, and by Lemma~\ref{lem:perp-linear} so is $U=W^{\perp'}$.  Therefore $a\mid\dimq U=n$, and Definition~\ref{def:expansion} shows that $\cC$ is, up to equivalence, the $\fq$-expansion of an $[n/a,k]_{q^m/q^a}$ code.  Since expansion multiplies all rank weights by $a$, every rank weight of $\cC$ is divisible by $a$.
\end{proof}

\begin{remark}\label{rem:square-case}
For two-dimensional codes, the square case is already covered by \cite[Proposition~IV.1 and Corollary~IV.4]{PolverinoEtAl2023}.  More precisely, let $\cC$ be a nondegenerate $[m,2]_{q^m/q}$ code and let $m-e$ be its second largest nonzero rank weight.  Then $e\mid m$, all rank weights of $\cC$ are divisible by $e$, and $\cC$ is equivalent to a code arising from $\mathbb F_{q^e}$.  Thus, in the square case, no divisor-gap assumption is required.
\end{remark}

\section{Codes with two or three nonzero weights}\label{sec:few-weight}

In this section we collect the consequences of the previous sections for codes with two or three nonzero weights.  By Proposition~\ref{prop:maximum-rank-weight} and Remark~\ref{rem:three-weight-maximum}, if a nondegenerate $[n,k]_{q^m/q}$ code has length $n\ge m$, then $m$ is one of its weights.  For such codes, Lemmas~\ref{lem:code-dichotomy} and~\ref{lem: equiv-cond-direct-sum} and Corollary~\ref{cor:direct-sum-weights} give a dichotomy:
\begin{enumerate}[label=\textup{(\Roman*)}]
  \item either $\cC_{<m}$ contains no codeword of rank $m$; equivalently, $h(\cC)=k-\lfloor n/m\rfloor$, and $\cC$ is the direct sum of a simplex code and a code of length smaller than $m$ (Section~\ref{subsec:direct-sum});
  \item or $\cC_{<m}$ contains codewords of rank $m$; equivalently, $n\ge(k-h(\cC)+1)m$, and the results of Section~\ref{subsec:not-direct-sum} apply.
\end{enumerate}
Two-dimensional codes of length $n\le m$ have no codeword of rank $m$ unless $n=m$; they are treated in Section~\ref{sec:antipodal}.  We first need two auxiliary results.

Lemma~\ref{lem:one-weight-short} can  be deduced from the classification of constant rank-weight codes in \cite{Randrianarisoa2020}.  

\begin{lemma}\label{lem:one-weight-short}
Let $\cC$ be a nondegenerate $[n,k]_{q^m/q}$ code with $n\le m-1$ whose nonzero codewords all have the same rank.  Then $k=1$.
\end{lemma}

The next lemma shows that, for linear sets with a single positive point weight, no divisor-gap condition is needed in Theorem~\ref{thm:gcd-rigidity}.

\begin{lemma}\label{lem:one-weight-linear-set}
Let $q\ge m$, let $U\le_{\fq}V$, and put $S=\langle U\rangle_{\fqm}$, $h=\dimqm S\ge2$ and $N=\dimq U\le(h-1)m$.  Assume that all points of $L_U$ have the same weight $a\ge2$.  Then $a\mid m$, $a\mid N$, $U$ is $\mathbb F_{q^a}$-linear and
\[
  \dim_{\mathbb F_{q^a}}(U\cap P)=1\qquad\text{for every } P\in L_U,
\]
that is, $U$ is $\mathbb F_{q^a}$-scattered.  Moreover, every divisor $t$ as in Theorem~\ref{thm:larger-field-representation} equals $a$.
\end{lemma}
\begin{proof}
    All the assertions except the last one are \cite[Corollary~III.3]{ZulloEtAl2024}, applied in $S$ instead of $V$.  Let $t$ and $W=\langle U\rangle_{\mathbb F_{q^t}}$ be as in Theorem~\ref{thm:larger-field-representation}, so that $L_W=L_U$, $\dimq W=N+t-a$ and $a\le t\le N-a$.  For $P\in L_U$, the quotient $(W\cap P)/(U\cap P)$ embeds into $W/U$, so $a\le w_W(P)\le a+(t-a)=t$; as $w_W(P)$ is a positive multiple of $t$, we get $w_W(P)=t$.  Applying \eqref{eq:vector-counting} to $U$ and to $W$,
\[
  \frac{q^N-1}{q^a-1}=|L_U|=|L_W|=\frac{q^{N+t-a}-1}{q^t-1},
\]
that is, $(q^{N-a}-1)(q^t-q^a)=0$.  Since $N-a\ge t\ge1$, we obtain $t=a$.
\end{proof}

\subsection{Two-weight codes}\label{subsec:two-weight}

The study of two-weight rank-metric codes was initiated by Pratihar and
Randrianarisoa~\cite{PratiharRandrianarisoa2024}, who considered \emph{antipodal}
two-weight codes. They proved that such codes have dimension two
\cite[Theorem~2]{PratiharRandrianarisoa2024} and that $n-d\mid n$, where $d<n$ is
the other nonzero weight \cite[Corollary~1]{PratiharRandrianarisoa2024}. Since every
nondegenerate code contains a codeword of rank $\min\{m,n\}$
\cite[Proposition~3.11]{AlfaranoEtAl2022}, the nonzero weights of a nondegenerate
two-weight code are $d$ and $\min\{m,n\}$. Hence two-weight codes of length $n\le m$
are antipodal and two-dimensional, so the restriction $k=2$ in
Theorem~\ref{thm:two-weight-short} below is no loss of generality. Two-weight codes
of length $n\ge m$ have nonzero weights $d<m$ and $m$; they are treated in
Theorem~\ref{thm:two-weight}.

For $q\ge m$, both cases were studied in~\cite{ZulloEtAl2024}. There it is shown
that antipodal two-weight codes are $\mathbb F_{q^{n-d}}$-scattered codes
\cite[Theorem~IV.8]{ZulloEtAl2024}, and that, for $n>m$, the system of a two-weight
code is the dual of an $\mathbb F_{q^{m-d}}$-scattered subspace
\cite[Corollary~IV.10 and Theorem~IV.11]{ZulloEtAl2024}.
Theorems~\ref{thm:two-weight-short} and~\ref{thm:two-weight}(ii) recover these
results, and remove the assumption on $q$ when $m\mid n$.

They also show that the second statement of~\cite{ZulloEtAl2024} requires an
additional hypothesis. The counting argument in the proofs of
\cite[Corollary~IV.10 and Theorem~IV.11]{ZulloEtAl2024} tacitly assumes
$h(\cC)\ge2$. When $h(\cC)=1$, it yields $km-n=m-d$ rather than a contradiction,
and the conclusion $m-d\mid m$ may fail. For instance, for every $q$ and every
$k\ge2$, the direct sum of a $[3(k-1),k-1]_{q^3/q}$ simplex code and a
nondegenerate $[1,1]_{q^3/q}$ code is a nondegenerate two-weight code with weights
$1$ and $3$, although $3-1\nmid3$. Theorem~\ref{thm:two-weight}(i) shows that
these direct sums are the only exceptions. Every two-weight code with $h(\cC)=1$ is
the direct sum of a $[(k-1)m,k-1]_{q^m/q}$ simplex code and a nondegenerate
$[d,1]_{q^m/q}$ code. Thus the classification of~\cite{ZulloEtAl2024} holds
precisely for the two-weight codes with $h(\cC)\ge2$, that is, for those whose
associated system contains no $\fqm$-hyperplane.

\begin{theorem}\label{thm:two-weight}
Let $\cC$ be a nondegenerate $[n,k]_{q^m/q}$ code whose nonzero weights are $d<m$ and $m$, let $U$ be an associated system, and put $a=m-d$ and $N=km-n$.  Then $a\mid N$,
\[
  A_d(\cC)=(q^m-1)\frac{q^N-1}{q^a-1},
  \qquad
  A_m(\cC)=q^{km}-1-A_d(\cC),
\]
and exactly one of the following holds.
\begin{enumerate}[label=\textup{(\roman*)}]
  \item $h(\cC)=1$.  Then $k\ge2$, $n=(k-1)m+d$ and, up to equivalence, $\cC=\tilde\cC\times\cC'$, where $\tilde\cC$ is a $[(k-1)m,k-1]_{q^m/q}$ simplex code and $\cC'$ is a nondegenerate $[d,1]_{q^m/q}$ code.  Conversely, for every $k\ge2$ and $1\le d\le m-1$, every such code is a nondegenerate two-weight code with weights $d$ and $m$ and with $h(\cC)=1$.
  \item $h(\cC)\ge2$.  Then $\cC_{<m}$ contains codewords of rank $m$, that is, $n\ge(k-h(\cC)+1)m$.  If moreover $q\ge m$ or $m\mid n$ or $d=n-1$, then
  \[
    a\mid m,\qquad a\mid n,\qquad a\mid d,
  \]
  $U$ and $U^{\perp'}$ are $\mathbb F_{q^a}$-linear, $U^{\perp'}$ is $\mathbb F_{q^a}$-scattered, and $\cC$ is, up to equivalence, the $\fq$-expansion of an $[n/a,k]_{q^m/q^a}$ code with nonzero weights $d/a$ and $m/a$.
\end{enumerate}
\end{theorem}

\begin{proof}
% Let $W = \Uperpp$. Then $\dimq W = N$ and since every point $P \in L_W\setminus\{0\}$ has weight $\wt(P) = m-d = a$, we have that $\lvert L_W \rvert = \dfrac{q^N - 1}{q^a - 1}$.
% So $q^a-1 \vert q^N-1$ and $a\vert N.$ Furthermore, $$A_d(\cC)= (q^m-1) \lvert L_W \rvert = (q^m-1)\frac{q^N-1}{q^a-1},$$ and $A_m(\cC)=q^{km}-1-A_d(\cC)$. \\
% (i) Let $h(\cC) = 1$. Then $\cC_{<m} = \{ c\in \cC : \wt(c) = d \}$ is generated by one codeword $c \in \cC$ of rank weight $\wt(c)< m$. 
% By Lemma \ref{lem: equiv-cond-direct-sum}, we get that $1 = k-\lfloor n/m\rfloor$. Since $\cC$ contains a codeword of weight $m$, we have $n \ge m$, so $k \ge 2.$
% By Corollary \ref{cor:direct-sum-weights}, we get up to equivalence, that $$\cC=\tilde\cC\times\cC',$$ where $\tilde\cC$ is a $[(k-1)m,k-1]_{q^m/q}$ simplex code and $\cC'$ is a nondegenerate $[\varepsilon,1]_{q^m/q}$ code with $\varepsilon=n-(k-1)m = d$.  \\
% (ii) Let $h(\cC) \ge 2$. Then $\dimqm(\cC_{<m}) \ge 2$, and since the only possible rank weights in $\cC$ are $d$ and $m$, we have that $\cC_{<m}$ contains codewords of rank $m$. By Lemma \ref{lem: equiv-cond-direct-sum} we get that $n \ge (k-h(\cC) + 1 )m.$ If moreover $q \ge m$ or  $m\mid n$ then the result follows from Corollary
% \ref{cor:code-rigidity} resp. Theorem \ref{thm:length-multiple-m}.

% \medskip

Let $W = \Uperpp$. Then $\dimq W = N$ and since, by \eqref{eq:dual-point-code-weight}, every point $P \in L_W$ has weight $\wt_W(P) = m-d = a$, we have, by \eqref{eq:vector-counting}, that $\lvert L_W \rvert = \dfrac{q^N - 1}{q^a - 1}$.
So $q^a-1 \vert q^N-1$ and $a\vert N.$ Furthermore, $$A_d(\cC)= (q^m-1) \lvert L_W \rvert = (q^m-1)\frac{q^N-1}{q^a-1},$$ and $A_m(\cC)=q^{km}-1-A_d(\cC)$. Since the weight $d$ occurs, $h(\cC)\ge1$, so exactly one of (i) and (ii) holds. \\
(i) Let $h(\cC) = 1$. Then $\cC_{<m} = \{ c\in \cC : \wt(c) = d \}\cup\{0\}$ is generated by one codeword $c \in \cC$ of rank weight $d=\wt(c)< m$, so it contains no codeword of rank $m$. 
By Lemma \ref{lem: equiv-cond-direct-sum}, (a)$\Rightarrow$(c), we get that $1 = k-\lfloor n/m\rfloor$. Since $\cC$ contains a codeword of weight $m$, we have $n \ge m$, so $k \ge 2.$
By Corollary \ref{cor:direct-sum-weights}, we get up to equivalence, that $$\cC=\tilde\cC\times\cC',$$ where $\tilde\cC$ is a $[(k-1)m,k-1]_{q^m/q}$ simplex code and $\cC'$ is a nondegenerate $[\varepsilon,1]_{q^m/q}$ code with $\varepsilon=n-(k-1)m = d$; indeed, the entries of a generator of $\cC'$ are $\fq$-linearly independent, so every nonzero codeword of $\cC'$ has rank $\varepsilon$, and by Corollary~\ref{cor:direct-sum-weights} these are exactly the codewords of $\cC$ of rank smaller than $m$, i.e.\ of rank $d$. 
Conversely, let $k\ge2$, $1\le d\le m-1$ and $\cC=\tilde\cC\times\cC'$ as above. The system of $\cC$ is $\fqm^{k-1}\times U'$, where $U'\le\fqm$ is the $d$-dimensional $\fq$-span of the entries of a generator of $\cC'$, so $\cC$ is nondegenerate. A codeword $(\tilde c,c')$ has rank $m$ if $\tilde c\ne0$, by Remark~\ref{rem:simplex-expansion}, and rank $\wt(c')=d$ if $\tilde c=0\ne c'$; both cases occur, and $\cC_{<m}=\{0\}\times\cC'$ is one-dimensional, so $h(\cC)=1$.\\
(ii) Let $h(\cC) \ge 2$. Then $\dimqm(\cC_{<m}) \ge 2$, and since the only possible rank weights in $\cC$ are $d$ and $m$, we have that $\cC_{<m}$ contains codewords of rank $m$. Indeed, put $h=h(\cC)$ and $S=\langle W\rangle_{\fqm}$, so that $\dimqm S=h$ (see the proof of Lemma~\ref{lem:code-dichotomy}). If $\cC_{<m}$ contained no codeword of rank $m$, then, by Lemma~\ref{lem: equiv-cond-direct-sum}, $L_W=\PG(S,\fqm)$ and all its points have weight $a$, so \eqref{eq:vector-counting} would give $q^N-1=\frac{q^{hm}-1}{q^m-1}(q^a-1)$, that is,
\[
q^N=q^a\Bigl(1+q^{m-a}\,\tfrac{q^{(h-1)m}-1}{q^m-1}\,(q^a-1)\Bigr).
\]
Since $1\le a\le m-1$ and $h\ge2$, the factor in brackets is larger than one and congruent to $1$ modulo $q$, a contradiction. By Lemma \ref{lem: equiv-cond-direct-sum} we get that $n \ge (k-h(\cC) + 1 )m$, that is, $N\le(h-1)m$. We may assume $a\ge2$, the case $a=1$ being trivial. If $q\ge m$, Lemma~\ref{lem:one-weight-linear-set}, applied to $W$ in $S$, shows that $a\mid m$ and that $W$ is an $\mathbb F_{q^a}$-scattered $\mathbb F_{q^a}$-subspace; if $m\mid n$, the same follows from Proposition~\ref{prop:multiple-m}, since $N=km-n$ is a multiple of $m$ and all points of $L_W$ have weight $a$. By Lemma~\ref{lem:perp-linear}, $U=W^{\perp'}$ is $\mathbb F_{q^a}$-linear; hence $a\mid n$, $a\mid d=m-a$ and, by Definition~\ref{def:expansion}, $\cC$ is, up to equivalence, the $\fq$-expansion of an $[n/a,k]_{q^m/q^a}$ code with nonzero weights $d/a$ and $m/a$.

\end{proof}

For two-dimensional codes of length at most $m$, the same argument applied to the system itself gives the following result.

\begin{theorem}\label{thm:two-weight-short}
Let $\cC$ be a nondegenerate $[n,2]_{q^m/q}$ code, with $n\le m$, whose nonzero weights are $d<n$ and $n$, let $U$ be an associated system and put $a=n-d$.  Then $a\mid n$ and
\[
  A_d(\cC)=(q^m-1)\frac{q^n-1}{q^a-1},
  \qquad
  A_n(\cC)=q^{2m}-1-A_d(\cC).
\]
If moreover $q\ge m$ or $n=m$ or $d=n-1$, then $a\mid m$, $a\mid d$, $U$ is an $\mathbb F_{q^a}$-scattered $\mathbb F_{q^a}$-subspace of $\fqm^2$, and $\cC$ is, up to equivalence, the $\fq$-expansion of an $[n/a,2]_{q^m/q^a}$ code with nonzero weights $d/a$ and $n/a$.  In particular, if $q\ge m$ and $a\ge2$, every divisor $t$ as in Theorem~\ref{thm:antipodal-puncturing} equals $a$, so that no coordinate has to be deleted.
\end{theorem}

\begin{proof}
% \textcolor{olive}{
%     Since $\dimq(U) = n$ and $\wt(c) = n - \dimq (U \cap x^\perp)$ for any codeword $c = xG$, we get that the hyperplane-intersection spectrum of the associated system U is given by $\{n-n,n-d\}$ or $\{0,a\}.$ 
%     Hence, for every point $P$ in $L_U$, we have $\wt(P) = \dimq (U\cap P) = a.$
%    % Let $N_a(U)$ be the set of hyperplanes in $\PG(1,q^m)$, with $\dimq(U\cap H) = a$. 
%    Hence $a\mid n$ and the number $A_d(\cC)$ of codewords of rank weight
% $d$ satisfies
% \[
% A_d(\cC)
% =
% (q^m-1)|L_U| = (q^m-1)\frac{q^n-1}{q^a-1}.
% \]
% If $n=m$, then we use the previous theorem \ref{thm:two-weight}.
% If $q\ge m$, we have the result in \cite[Theorem IV.8]{ZulloEtAl2024} showing that $a\mid m$, $a\mid d$ and $U$ is a $\mathbb{F}_{q^a}$-scattered $\mathbb{F}_{q^a}$-subspace of $\fqm^2$. In particular, we have that $\cC$, up to equivalence, is the $\fq$-expansion of an $[n/a,2]_{q^m/q^a}$ code with nonzero weights $d/a$ and $n/a$. 
% }

% \medskip

    Since $\dimq(U) = n$ and $\wt(c) = n - \dimq (U \cap x^\perp)$ for any codeword $c = xG$, we get that the hyperplane-intersection spectrum of the associated system U is given by $\{n-n,n-d\}$ or $\{0,a\}.$ 
    Hence, for every point $P$ in $L_U$, we have $\wt(P) = \dimq (U\cap P) = a.$
   % Let $N_a(U)$ be the set of hyperplanes in $\PG(1,q^m)$, with $\dimq(U\cap H) = a$. 
   Hence $q^n-1=|L_U|(q^a-1)$ by \eqref{eq:vector-counting}, so $a\mid n$ and the number $A_d(\cC)$ of codewords of rank weight
$d$ satisfies
\[
A_d(\cC)
=
(q^m-1)|L_U| = (q^m-1)\frac{q^n-1}{q^a-1}.
\]
Moreover, $A_n(\cC)=q^{2m}-1-A_d(\cC)$ by \eqref{eq:antipodal-counting}.
If $n=m$, then we use the previous Theorem \ref{thm:two-weight}, with $a=m-d=n-d$: here $h(\cC)=2$, since $h(\cC)=1$ would give $n=m+d>m$ by Theorem~\ref{thm:two-weight}(i), so case~(ii) applies with $m\mid n$ and shows that $a\mid m$, $a\mid d$ and $U$ is $\mathbb F_{q^a}$-linear; as every point of $L_U$ has weight $a$, $U$ is $\mathbb F_{q^a}$-scattered.
If $q\ge m$, we have the result in \cite[Theorem IV.8]{ZulloEtAl2024} (whose hypothesis $n=\ell(n-d)$ is $a\mid n$) showing that $a\mid m$, $a\mid d$ and $U$ is a $\mathbb{F}_{q^a}$-scattered $\mathbb{F}_{q^a}$-subspace of $\fqm^2$. In particular, we have that $\cC$, up to equivalence, is the $\fq$-expansion of an $[n/a,2]_{q^m/q^a}$ code with nonzero weights $d/a$ and $n/a$. Finally, if $q\ge m$ and $a\ge2$, the divisor $t$ in Theorem~\ref{thm:antipodal-puncturing} is one given by Theorem~\ref{thm:larger-field-representation} applied to $U$ (step~1 of its proof), so $t=a$ by the last assertion of Lemma~\ref{lem:one-weight-linear-set}, and no coordinate has to be deleted.
\end{proof}

\begin{remark}\label{rem:two-weight-scattered}
Theorems~\ref{thm:two-weight} and~\ref{thm:two-weight-short} show that no divisor-gap condition is needed for two-weight codes.  Let $a$ be the difference between the two nonzero weights.  If $q\ge m$, every nondegenerate two-weight code of length $n\ge m$ with $h(\cC)\ge2$, as well as every nondegenerate two-weight $[n,2]_{q^m/q}$ code with $n\le m$, is, up to equivalence, an $\fq$-expansion of a code over $\mathbb F_{q^a}$ whose associated $\mathbb F_{q^a}$-linear set is $\mathbb F_{q^a}$-scattered \textup($U^{\perp'}$ in Theorem~\ref{thm:two-weight}, $U$ in Theorem~\ref{thm:two-weight-short}\textup); the remaining two-weight codes of length $n\ge m$ are the direct sums of Theorem~\ref{thm:two-weight}(i).  Compare with the description of two-weight rank-metric codes through scattered subspaces in \cite{ZulloEtAl2024}.  In contrast, Family~$\mathcal F_4$ with $t\geq 3$ (Example~\ref{ex:F4}; for instance $m=6, t=3,\ell=2$, which gives the weights $2,3,5$) shows that the analogous statement fails for three weights: a three-weight $[n,2]_{q^m/q}$ code with $n\le m$ and first gap at least two need not be an expansion.
\end{remark}

\subsection{Three-weight codes}\label{subsec:three-weight}

For three-weight codes of length $n\ge m$, the weights are $d<w<m$ together with $m$.  We treat the two cases of the dichotomy separately.

\begin{theorem}\label{thm:three-weight-direct}
Let $\cC$ be a nondegenerate $[n,k]_{q^m/q}$ code, with $n\ge m$, whose nonzero weights are $d<w<m$ and $m$, and assume that $\cC_{<m}$ contains no codeword of rank $m$.  Put $h=h(\cC)$ and $\varepsilon=n-(k-h)m$.  Then $2\le h\le k-1$, $0<\varepsilon\le m-1$ and, up to equivalence, $\cC=\tilde\cC\times\cC'$, where $\tilde\cC$ is a $[(k-h)m,k-h]_{q^m/q}$ simplex code and $\cC'$ is a nondegenerate $[\varepsilon,h]_{q^m/q}$ code with nonzero weights $d<w\le\varepsilon$.  The weight distribution of $\cC$ is determined by $q,m,n,k,d,w$: putting $\theta=\frac{q^{hm}-1}{q^m-1}$ and $N=km-n=hm-\varepsilon$,
\begin{align*}
  A_d(\cC)&=\frac{q^m-1}{q^{m-d}-q^{m-w}}\Bigl[q^N-1-\theta\bigl(q^{m-w}-1\bigr)\Bigr],\\
  A_w(\cC)&=\frac{q^m-1}{q^{m-d}-q^{m-w}}\Bigl[\theta\bigl(q^{m-d}-1\bigr)-\bigl(q^N-1\bigr)\Bigr],\\
  A_m(\cC)&=q^{km}-q^{hm}.
\end{align*}
If $h=2$, then $w=\varepsilon=n-(k-2)m$; if moreover $q\ge m$, then $\varepsilon-d$ divides $m$, $\varepsilon$ and $d$, and $\cC'$ is, up to equivalence, the $\fq$-expansion of an $\bigl[\frac{\varepsilon}{\varepsilon-d},2\bigr]_{q^m/q^{\varepsilon-d}}$ code with nonzero weights $\frac d{\varepsilon-d}$ and $\frac{\varepsilon}{\varepsilon-d}$.
\end{theorem}

\begin{proof}
By Lemma~\ref{lem: equiv-cond-direct-sum}(c), $h=k-\lfloor n/m\rfloor\le k-1$, as $n\ge m$, and by Corollary~\ref{cor:direct-sum-weights}, $\cC=\tilde\cC\times\cC'$ with $\cC'$ nondegenerate of length $\varepsilon\le m-1$, whose nonzero weights are exactly $d$ and $w$.  In particular $\varepsilon\ge w>0$.  If $h=1$, all nonzero codewords of $\cC'$ would have the same rank; hence $h\ge2$.

Let $W=U^{\perp'}$ for an associated system $U$.  By Lemma~\ref{lem: equiv-cond-direct-sum}(d), $L_W$ is the whole projective space $\PG(\langle W\rangle_{\fqm},\fqm)\cong\PG(h-1,q^m)$, whose points have weight $m-d$ or $m-w$ by \eqref{eq:dual-point-code-weight}.  If $x$ and $y$ denote the numbers of points of weight $m-d$ and $m-w$, respectively, then $x+y=\theta$ and, by \eqref{eq:vector-counting}, $x(q^{m-d}-1)+y(q^{m-w}-1)=q^N-1$.  Solving this system and using $A_d(\cC)=(q^m-1)x$, $A_w(\cC)=(q^m-1)y$ gives the formulas; $A_m(\cC)$ is given by Corollary~\ref{cor:direct-sum-weights}.

If $h=2$, then $\cC'$ is a nondegenerate $[\varepsilon,2]_{q^m/q}$ code with $\varepsilon\le m$, hence it is antipodal (see the beginning of Section~\ref{sec:antipodal}), so $\varepsilon$ is one of its weights and $w=\varepsilon$.  The last assertion follows from Theorem~\ref{thm:two-weight-short} applied to $\cC'$.
\end{proof}

When $\cC_{<m}$ contains codewords of rank $m$, the results of Section~\ref{subsec:not-direct-sum} specialize as follows.

\begin{theorem}\label{thm:three-weight-not-direct}
Let $\cC$ be a nondegenerate $[n,k]_{q^m/q}$ code whose nonzero weights are $d<w<m$ and $m$, and assume that $\cC_{<m}$ contains codewords of rank $m$, that is, $n\ge(k-h(\cC)+1)m$.  Put $N=km-n$.  Then $h(\cC)\ge2$, and the following hold.
\begin{enumerate}[label=\textup{(\roman*)}]
  \item Assume $q\ge m$ and $m-w\ge2$.  Then there exists a proper divisor $t$ of $m$ with
  \[
    m-w\le t\le N-(m-w),
    \qquad
    t\mid n-w,
    \qquad
    d\bmod t\le t-(m-w),
  \]
  and $w\bmod t=t-(m-w)$ if $t>m-w$.  Deleting $t-(m-w)$ suitable coordinates of $\cC$ yields the $\fq$-expansion of a code over $\mathbb F_{q^t}$, and a codeword of rank $d$, $w$ or $m$ becomes a codeword of rank $t\lfloor d/t\rfloor$, $m-t$ or $m$, respectively.
  \item Assume $q\ge m$ and $m-w\ge2$, put $e=\gcd(m-w,m-d)$, and assume $\cD_m(e,N)=\varnothing$.  Then $e=m-w$; consequently
  \[
    m-w\mid m-d,\qquad
    m-w\mid m,\qquad
    m-w\mid n,
  \]
  and $\cC$ is, up to equivalence, the $\fq$-expansion of an $\left[\frac{n}{m-w},k\right]_{q^m/q^{m-w}}$ code with nonzero weights $\frac d{m-w}$, $\frac w{m-w}$ and $\frac m{m-w}$.
  \item If $m\mid n$, then the conclusion of~(ii) holds without any assumption on $q$, on $m-w$ or on the divisors of $m$.
\end{enumerate}
\end{theorem}

\begin{proof}
Since $\cC$ has two distinct nonmaximum weights, $h(\cC)\ge2$ by Remark~\ref{rem:code-hypotheses}.  Under the assumptions of~(i) and~(ii), \ref{H1}--\ref{H4} of Section 5.2 hold with $r=2$, $d_1=d$, $d_2=w$ and $a=m-w$.

(i) Theorem~\ref{thm:larger-field-puncturing} and Proposition~\ref{prop:weight-transformation}(ii) give the divisor $t$ and all the stated properties except properness; $t\ne m$ by Remark~\ref{rem:refined-gap}(i), applied to $U^{\perp'}$.  By Theorem~\ref{thm:larger-field-puncturing}(ii), codewords of rank $d$, $w$ and $m$ become codewords of rank $t\lfloor d/t\rfloor$, $t\lfloor w/t\rfloor$ and $m$, and $t\lfloor w/t\rfloor=m-t$ because $w=m-a$ with $1\le a\le t$ and $t\mid m$.

(ii) This is Corollary~\ref{cor:three-weight-rigidity}.

(iii) This is Theorem~\ref{thm:length-multiple-m} with $d_r=w$: it gives $m-w\mid m$, $m-w\mid n$ and the expansion; as all weights are divisible by $m-w$, also $m-w\mid d$.
\end{proof}

\begin{remark}\label{rem:three-weight-summary}
Theorems~\ref{thm:three-weight-direct} and~\ref{thm:three-weight-not-direct} cover all nondegenerate three-weight codes of length $n\ge m$.  In case~(I), the problem is reduced to two-weight codes of length smaller than $m$, and for $h=2$ such codes are completely described by Theorem~\ref{thm:two-weight-short} when $q\ge m$.  In case~(II), without the divisor-gap condition, some coordinates may have to be deleted before obtaining an expansion; see Example~\ref{ex:weight-collapse}.  For three-weight codes of length $n<m$, which have no codeword of rank $m$, we only considered the two-dimensional case, discussed below.
For $n=m$ and $k=2$, Theorem~\ref{thm:three-weight-not-direct}(iii) applies for every $q$;  indeed $h(\cC)=2$, since if $h(\cC)=1$, then $U^{\perp'}$, of $\fq$-dimension $m$, would be a point $P$, and $U=P^{\perp}$ would not span $\fqm^2$; hence $N=m\le(h(\cC)-1)m$, and $\cC_{<m}$ contains codewords of rank $m$ by Lemma~\ref{lem:code-dichotomy}.
\end{remark}

\begin{remark}
Corollary~\ref{cor:code-rigidity}, Theorem~\ref{thm:length-multiple-m}, Remark~\ref{rem:antipodal-consequences}(d) and Theorems~\ref{thm:two-weight} and~\ref{thm:two-weight-short} give positive parameter regimes for the problem studied in \cite{PolverinoEtAl2023}: when does an $e$-divisible $\mathbb F_{q^m}$-linear rank-metric code arise from a code over $\mathbb F_{q^e}$?  Our hypotheses are stated directly in terms of the weight spectrum and the divisors of $m$. 
{Already for two-dimensional codes of length $n\le m$, the families of Section~\ref{subsec:antipodal-examples} show that three-weight codes are much less rigid than two-weight codes.  By Theorem~\ref{thm:two-weight-short}, if $q\ge m$, a nondegenerate two-weight $[n,2]_{q^m/q}$ code with weights $d<n$ is completely determined, up to the choice of an $\mathbb F_{q^a}$-scattered subspace, by its parameters: the gap $a=n-d$ divides $m$, $n$ and $d$, the code is the $\fq$-expansion of a code over $\mathbb F_{q^a}$ (no coordinate has to be deleted), and its weight distribution depends only on $q$, $m$, $n$ and $a$.  For three weights none of these properties survives, and the constructions of Section~\ref{subsec:antipodal-examples}, which exist for every prime power $q$ (in particular for $q\ge m$), exhibit quite different behaviours (see also the table in Remark~\ref{rem:gap-behaviour}):
\begin{itemize}
  \item Family~$\mathcal F_1$ (Example~\ref{ex:F1}): weights $n-2,n-1,n$ for every $3\le n\le m$, with $g_1=g_2=1$; here the counting relations do \emph{not} determine the weight distribution from the weights, since $A_{n-2}(\cC)=(q^m-1)\nu$ depends on the number $\nu\ge1$ of points of weight two, which is not fixed by the construction.
  \item Family~$\mathcal F_2$ (Example~\ref{ex:F2}) and, more generally, $i$-clubs (Remark~\ref{rem:clubs-general}): weights $N-i,N-1,N$, with $g_1=1$ and $g_2=i-1$, for $2\le i\le N-1$ and $N\le m$; the second gap can be as large as $m-3$ (Family~$\mathcal F_2$ with $n=m$), the minimum weight can be as small as $2$, or even $1$ (clubs obtained from $\{(x,\Tr_{q^m/q}(x))\}$), and there is, up to $\fqm$-multiples, a single codeword of minimum weight, $A_{N-i}(\cC)=q^m-1$.  No divisibility among $N$, $m$ and the weights is forced.
  \item Family~$\mathcal F_3$ (Example~\ref{ex:F3}): weights $a,2a-1,2a$, for $m$ even and $2\le a\le m/2$, with $g_1=1$ and $g_2=a-1$; the second largest weight can be coprime with $m$ (weights $3,5,6$ for $m=6$, $a=3$).
  \item Family~$\mathcal F_4$ (Example~\ref{ex:F4}): weights $\ell t-t-1,\ \ell t-t,\ \ell t-1$, for $t\mid m$, $1<t<m$ and $2\le\ell\le m/t$, with $g_1=t-1$ and $g_2=1$; the first gap can be as large as $m/2-1$, and yet $\cC$ is not an expansion, since $\gcd(g_1,g_1+g_2)=1$: it is obtained by deleting one coordinate from an expansion over $\mathbb F_{q^t}$.  This is exactly the phenomenon described by Theorem~\ref{thm:antipodal-puncturing} with $t>a$, which Theorem~\ref{thm:two-weight-short} rules out for two-weight codes.
\end{itemize}
Thus, while a gap at least two in a two-weight code forces the code to come from a subfield, in a three-weight code a large gap is compatible both with the absence of any divisibility condition (second gap, Families~$\mathcal F_2$, $\mathcal F_3$ and clubs) and with a code that is only a punctured expansion (first gap, Family~$\mathcal F_4$).  Even when $q\ge m$, the only general restrictions known are those of Corollary~\ref{cor:antipodal-puncturing-three}, which come into play only when the first gap is at least two.}
\end{remark}

\section*{Declaration of Generative AI and AI-assisted Technologies in the Writing Process}

The authors used AI-assisted tools exclusively for language editing and improving the clarity and readability of the manuscript. All mathematical results, proofs, and scientific conclusions are the authors' own.
The authors have carefully reviewed the entire manuscript and take full responsibility for its content.

\section*{Acknowledgements}

The research of Olga Polverino and Ferdinando Zullo was partially supported by the
Italian National Group for Algebraic and Geometric Structures and their Applications
(GNSAGA - INdAM). This research was also partially supported by the 
Scientific Research Community of the Research Foundation Flanders (Belgium), {\em
Graphs, Association schemes and Geometries: structures, algorithms and computation
(W003324N)} and by the VUB-OZR project {\em Discrete Structures, and
their applications in Data Science” (OZR3637)}. This work was partially written during
the Opera 2026 conference in Bordeaux, and so we acknowledge the support from the
International Research Laboratory LYSM in partnership between CNRS and INdAM.

\end{document}